\documentclass{amsart}
\usepackage{graphicx}    
\usepackage{color}
\usepackage{xspace}
\usepackage{hyperref}
\usepackage{esint}
\usepackage{graphicx}
\usepackage{amsfonts,amsmath,amssymb}
\usepackage{mathrsfs,mathtools}
\usepackage{enumerate}
\DeclareMathAlphabet{\mathpzc}{OT1}{pzc}{m}{it}

\usepackage{stmaryrd}
\usepackage{bbm}

\newcommand{\srn}{{\tfrac12}}

\newcommand{\R}{\mathbb{R}}

\newcommand{\Y}{\mathpzc{Y}}

\newcommand{\N}{\mathpzc{N}}

\newcommand{\C}{\mathcal{C}}

\newcommand{\V}{\mathbb{V}}

\newcommand{\Xcal}{\mathcal{X}}

\newcommand{\calI}{\mathcal{I}}
\newcommand{\T}{\mathscr{T}}
\newcommand{\E}{\mathscr{E}}

\newcommand{\Ss}{\mathscr{S}}

\newcommand{\HL}{ \mbox{ \raisebox{7.2pt} {\tiny$\circ$} \kern-10.7pt} {H_L^1} }
\newcommand{\Wp}{ \mbox{ \raisebox{7.7pt} {\scriptsize$\circ$} \kern-10.1pt} {W^{1,p}} }
\newcommand{\Wpp}{ \mbox{ \raisebox{7.7pt} {\scriptsize$\circ$} \kern-10.1pt} {W^{1,p'}} }
\newcommand{\Sz}{ \mbox{ \raisebox{7.5pt} {\scriptsize$\circ$} \kern-10.1pt} {\Ss} }
\newcommand{\HLnew}{ \mbox{ \raisebox{7pt} {\scriptsize$\circ$} \kern-10.1pt}{H}^1_L }
\newcommand{\HLn}{{\mbox{\,\raisebox{4.7pt} {\tiny$\circ$} \kern-9.3pt}{H}^{1}_{L}  }}
\newcommand{\HLs}{{\mbox{\raisebox{8.7pt} {\scriptsize$\circ$} \kern-10.1pt}{H}^1_L  }}

\newcommand{\Tr}{\mathbb{T}}
\newcommand{\U}{\mathbb{U}}

\newcommand{\usf}{\mathsf{u}}

\newcommand{\ue}{\mathscr{U}}

\DeclareMathOperator*{\tr}{tr_\Omega}

\newcommand{\boxednumber}[1]{\expandafter\readdigit\the\numexpr#1\relax\relax}
\newcommand{\K}{\mathcal{K}}
\newcommand{\resto}{\mathsf{r}}
\newcommand{\Hsd}{\mathbb{H}^{-s}(\Omega)}

\newcommand{\M}{\mathpzc{M}}
\newcommand{\Laps}{(-\Delta)^s}

\newcommand{\DIV}{\textrm{div}}
\newcommand{\diff}{\, \mbox{\rm d}}
\newcommand{\ie}{i.e.,\@\xspace}
\newcommand{\cf}{cf.\@\xspace}
\newcommand{\Hs}{\mathbb{H}^s(\Omega)}
\newcommand{\Ws}{\mathbb{H}^{1-s}(\Omega)}

\DeclareMathOperator*{\Span}{span}

\newcommand{\Nin}{\,{\mbox{\,\raisebox{6.0pt} {\tiny$\circ$} \kern-10.9pt}\N }}

\newcommand{\osc}{{\textup{\textsf{osc}}}}

\author[R.H.~Nochetto]{Ricardo H.~Nochetto}
\address[R.H.~Nochetto]{Department of Mathematics and Institute for Physical Science and Technology,
University of Maryland, College Park, MD 20742, USA}
\email{rhn@math.umd.edu}
\thanks{RHN has been supported in part by NSF grants DMS-1109325 and DMS-1411808.}

\author[E.~Ot\'arola]{Enrique Ot\'arola}
\address[E.~Ot\'arola]{Departamento de Matem\'atica, Universidad T\'ecnica Federico Santa Mar\'ia, Valparaiso, Chile}
\email{enriqueotarola@gmail.com}
\thanks{EO has been supported in part by NSF grants DMS-1109325 and DMS-1411808 and by CONICYT through project Anillo ACT1106.}

\author[A.J.~Salgado]{Abner J.~Salgado}
\address[A.J.~Salgado]{Department of Mathematics, University of Tennessee, Knoxville, TN 37996, USA}
\email{asalgad1@utk.edu}
\thanks{AJS has been supported in part by NSF grant DMS-1418784.}

\numberwithin{equation}{section}
\numberwithin{figure}{section}
\numberwithin{table}{section}
\newtheorem{theorem}{Theorem}[section]
\newtheorem{corollary}[theorem]{Corollary}
\newtheorem{lemma}[theorem]{Lemma}

\theoremstyle{definition}
\newtheorem{definition}[theorem]{Definition}
\newtheorem{remark}[theorem]{Remark}

\title[Numerical Fractional Diffusion]{A PDE Approach to Numerical Fractional Diffusion}

\begin{document}

\begin{abstract}
Fractional diffusion has become a fundamental tool for the modeling of
multiscale and heterogeneous phenomena. However, due to its nonlocal
nature, its accurate numerical approximation is delicate. We
survey our research program on the design and analysis of
efficient solution techniques for problems involving fractional
powers of elliptic operators.
Starting from a localization PDE result for these operators, 
we develop local techniques for their solution: a priori and a
posteriori error analyses, adaptivity and multilevel methods. 
We show the flexibility of our approach by proposing and
analyzing local solution techniques for a space-time fractional parabolic equation.
\end{abstract}



\date{The final version of this overview appeared in the Proceedings of the ICIAM, 2015}

\maketitle

\section{Introduction}

Diffusion is the tendency of a substance to evenly spread into
surrounding space, and is one of the most common physical
processes. The classical models of diffusion lead to local and
thoroughly studied equations. However, in recent times, it has become
evident that many of the assumptions that lead to these models are not
always satisfactory or not even realistic in practice. Consequently, different models of diffusion have been proposed, fractional diffusion being one of them.

Capturing the essential behavior of fractional diffusion with the simplest and crudest models is of paramount importance in science and engineering. This allows for understanding of physical applications where long range or anomalous diffusion is considered \cite{Abe2005403}, complex phenomena in mechanics \cite{atanackovic2014fractional}, biophysics \cite{bio}, turbulence \cite{wow}, image processing \cite{GH:14}, nonlocal electrostatics \cite{ICH}, finance \cite{MR2064019}
and the control of devices \cite{AO,AOS}. To understand the behavior
of fractional diffusion, computational science is fundamental. It is
one of the pillars, together with theory and experiments, of
scientific inquiry. A carefully crafted computational model can
replace a very expensive or unrealizable experimental setting, and it
can give new insight into the theoretical developments of a specific
discipline. The analysis of such computational schemes is the realm of
numerical analysis, which offers a rigorous mathematical
description of the extent to which the computer's output approximates the process of interest. These crucial aspects of modern research are blended together in this paper, which deals with the design and analysis of efficient numerical techniques for problems involving fractional diffusion.

The mathematical structure of fractional diffusion is shared by a wide
class of \emph{nonlocal operators} \cite{DuGuLe:12,Lipton:14}. These operators have a
strong connection with real-world problems, since they constitute a
fundamental part of the modeling and simulation of complex
phenomena. It is evident that the particular type of mathematical
operator appearing in applications can widely vary and that a
unified analysis might be well beyond reach. A more modest, but
nevertheless quite ambitious, goal is to develop computational tools
and their analysis for a representative of a particular class.
We discuss below our contributions to fractional diffusion.

Exploiting the breakthrough by L.~Caffarelli and L.~Silvestre
\cite{CS:07}, we have made a decisive advance in numerical
fractional diffusion. Although the analysis of our method is
intricate, \emph{its implementation is done using standard
components of finite element analysis}
\cite{CNOS2,CNOS,NOS4,NOS,NOS3}. This is the main advantage of our
scheme, since alternative approaches require less
traditional techniques (such as special quadrature) to cope with the
mathematical difficulties inherent to fractional diffusion. Even in 1D
\cite{HO,MR2684296}, 
the integral formulation of fractional diffusion is notoriously
difficult from the numerical standpoint due to the presence of a
nonintegrable kernel; see \cite{Otarola} for a
discussion. In contrast, our PDE approach to fractional
diffusion, as well as the recent method of \cite{BoPa:15},
can handle multidimensions easily and efficiently -- a highly
desirable feature.
  
Here we are interested in the design and analysis of numerical
techniques to solve problems involving fractional powers of the
Dirichlet Laplace operator $\Laps$, $s \in (0,1)$, the so-called
fractional Laplacian. Let $\Omega$ be a bounded Lipschitz domain
of $\R^n$ ($n\ge1$), with boundary $\partial\Omega$. Given a function
$f$, we seek $u$ such that
\begin{equation}
\label{fl=f_bdddom}
    \Laps u = f \quad \text{in } \Omega.
\end{equation}
Our approach, however, is not particular to the fractional Laplacian and it
can be applied to any second order, symmetric and uniformly elliptic operator \cite[\S7]{NOS}.

One of the main difficulties to study \eqref{fl=f_bdddom} is
that the fractional Laplacian is a nonlocal operator
\cite{CT:10,CDDS:11,CS:11,CS:07,Landkof}. To localize it, Caffarelli
and Silvestre showed in \cite{CS:07} that any power $s<1$ of the
fractional Laplacian in $\R^n$ can be realized as an operator that
maps a Dirichlet boundary condition to a Neumann-type condition via an
extension problem on the upper half-space $\R^{n+1}_+$. This
result was adapted in \cite{ BCdPS:12,CDDS:11,ST:10} to bounded
domains $\Omega$, thus obtaining an extension  problem posed on the semi-infinite cylinder $\C := \Omega \times (0,\infty)$. This extension is the mixed boundary value problem:
\begin{equation}
\label{alpha_harm_intro}
  \DIV\left( y^\alpha \nabla \ue \right) = 0  \text{ in } \C, \quad
  \ue = 0  \text{ on } \partial_L \C, \quad
  \frac{ \partial \ue }{\partial \nu^\alpha} = d_s f  \text{ on } \Omega \times \{0\},
\end{equation}
where $\partial_L \C := \partial \Omega \times (0,\infty)$ is the lateral boundary of $\C$, and $d_s: = 2^{1-2s} \frac{\Gamma(1-s)}{\Gamma(s)}$ is a positive normalization constant that depends only on $s$; see \cite{CS:11,CS:07} for details. The parameter $\alpha$ is defined as $\alpha = 1-2s \in (-1,1)$, and the so-called conormal exterior derivative of $\ue$ at $\Omega \times \{ 0 \}$ is
\begin{equation}
\label{def:lf}
\frac{\partial \ue}{\partial \nu^\alpha} = -\lim_{y \rightarrow 0^+} y^\alpha \ue_y.
\end{equation}
We call $y$ the \emph{extended variable} and the dimension $n+1$ in $\R_+^{n+1}$ the \emph{extended dimension} of problem \eqref{alpha_harm_intro}. The limit in \eqref{def:lf} must be understood in the distributional sense; see \cite{CS:11,CS:07,CDDS:11}. As noted in \cite{BCdPS:12,CS:07,CDDS:11, ST:10}, the fractional Laplacian and the Dirichlet-to-Neumann operator of problem \eqref{alpha_harm_intro} are related by
\begin{equation}
\label{eq:identity}
  d_s \Laps u = \frac{\partial \ue}{\partial \nu^\alpha } \quad \text{in } \Omega.
\end{equation}

We propose, and analyze in Section~\ref{sec:PDE_approach}, the following strategy to solve \eqref{fl=f_bdddom}: given $f$ we solve \eqref{alpha_harm_intro}, thus obtaining a function $\ue$; setting $u: x' \in \Omega \mapsto u(x') = \ue(x',0) \in \R$, we obtain  the solution of \eqref{fl=f_bdddom}.
The main advantage of this algorithm is that we are solving the local
problem \eqref{alpha_harm_intro} instead of dealing with the nonlocal
operator $\Laps$ of \eqref{fl=f_bdddom}. However, this comes at the
expense of incorporating one more dimension to the problem, thus
raising the question of computational efficiency. This motivates
the use of anisotropic meshes to compensate for the singular
behavior of $\ue(\cdot,y)$ as $y\to0$, and well as the development
of a posteriori error estimators and multilevel methods. These are
reviewed in Sections \ref{sec:a_posteriori} and \ref{sec:multilevel},
respectively. To show the flexibility of our approach, in
Section~\ref{sec:time_dependent} we consider a parabolic equation with fractional diffusion and fractional time derivative.


Throughout this work $\Omega$ is an open, bounded and connected domain
of $\R^n$, $n\geq1$, with polyhedral boundary
$\partial\Omega$. Given the semi-infinite cylinder $\C= \Omega\times (0,\infty)$
and $\Y>0$, the truncated cylinder with base $\Omega$ and height $\Y$ is defined by
$ \C_\Y := \Omega \times (0,\Y)$ with lateral boundary
$\partial_L\C_\Y := \partial \Omega \times (0,\Y)$.

We will be dealing with objects defined in $\R^{n+1}$ and it will be convenient to distinguish the extended dimension. A vector $x\in \R^{n+1}$, will be denoted by
\[
  x =  (x_1,\ldots,x_n, x_{n+1}) = (x', x_{n+1}) = (x',y),
\]
with $x_i \in \R$ for $i=1,\ldots,{n+1}$, $x' \in \R^n$ and $y\in\R$.

By $a \lesssim b$ we mean $a \leq Cb$, with a constant $C$ that neither 
depends on $a,b$ or the discretization parameters. The value of $C$ might change at each occurrence. 

\section{A PDE approach: formulation and FEM}
\label{sec:PDE_approach}
For a bounded domain there are several ways, not necessarily
equivalent, to define the fractional Laplacian; see \cite{BSV:14,NOS}
for a discussion. As in \cite{NOS} we adopt the definition based on spectral theory \cite{BS}.
The operator $(-\Delta)^{-1}: L^2(\Omega)\to L^2(\Omega)$ is compact, symmetric and positive, whence its spectrum $\{ \lambda_k^{-1} \}_{k\in \mathbb N}$ is discrete, real, positive and accumulates at zero. Moreover, there exists $\{ \varphi_k \}_{k\in \mathbb N}$, which is an orthonormal basis of $L^2(\Omega)$ and satisfies $- \Delta \varphi_k = \lambda_k \varphi_k$ in $\Omega$ and $\varphi_k = 0$ on $\partial\Omega$.
Fractional powers of the Dirichlet Laplacian can be defined for $w \in C_0^{\infty}(\Omega)$ by
\begin{equation}
  \label{def:second_frac}
  (-\Delta)^s w  = \sum_{k=1}^\infty \lambda_k^{s} w_k \varphi_k,
\end{equation} 
where $w_k = \int_{\Omega} w \varphi_k $. By density $(-\Delta)^s$ can be extended 
to the space
\begin{equation}
\label{def:Hs}
  \Hs = \left\{ w = \sum_{k=1}^\infty w_k \varphi_k: 
  \sum_{k=1}^{\infty} \lambda_k^s w_k^2 < \infty \right\}
\end{equation}
For $ s \in (0,1)$ we denote by $\Hsd$ the dual of $\Hs$.

\subsection{The Caffarelli-Silvestre extension problem}
\label{sub:CaffarelliSilvestre}
The Caffarelli-Silvestre extension leads to a \emph{nonuniformly} elliptic
equation \cite{CS:07,BCdPS:12,CT:10, CDDS:11}. We consider weighted Sobolev spaces with the weight $|y|^{\alpha}$, $\alpha \in (-1,1)$. If $D \subset \R^{n+1}$, $L^2(|y|^\alpha,D)$ is the space of measurable functions on $D$ such that
\[
\| w \|_{L^2(|y|^{\alpha},D)}^2 = \int_{D}|y|^{\alpha} w^2 < \infty.
\]
Similarly we define the weighted Sobolev space
\[
H^1(|y|^{\alpha},D) =
  \left\{ w \in L^2(|y|^{\alpha},D): | \nabla w | \in L^2(|y|^{\alpha},D) \right\},
\]
where $\nabla w$ is the distributional gradient of $w$. We equip $H^1(|y|^{\alpha},D)$ with the norm
\begin{equation}
\label{wH1norm}
\| w \|_{H^1(|y|^{\alpha},D)} =
\left(  \| w \|^2_{L^2(|y|^{\alpha},D)} + \| \nabla w \|^2_{L^2(|y|^{\alpha},D)} \right)^{1/2}.
\end{equation}
Since $\alpha \in (-1,1)$ then $|y|^\alpha$ belongs to the Muckenhoupt class $A_2(\R^{n+1})$ \cite{Javier,FKS:82,GU,Muckenhoupt,Turesson}. This, in particular, implies that $H^1(|y|^{\alpha},D)$ with norm \eqref{wH1norm} is Hilbert and $C^{\infty}(D) \cap H^1(|y|^{\alpha},D)$ is dense in $H^1(|y|^{\alpha},D)$ (cf.~\cite[Proposition 2.1.2, Corollary 2.1.6]{Turesson}, \cite{KO84} and \cite[Theorem~1]{GU}). We recall the definition of Muckenhoupt classes.

\begin{definition}[Muckenhoupt class $A_2$, \cite{Muckenhoupt,Turesson}]
 \label{def:Muckenhoupt}
Let $\omega$ be a weight and $N \geq 1$. We say that $\omega \in A_2(\R^N)$ if
\begin{equation}
  \label{A_pclass}
  C_{2,\omega} = \sup_{B} \left( \fint_{B} \omega \right)
            \left( \fint_{B} \omega^{-1} \right) < \infty,
\end{equation}
where the supremum is taken over all balls $B$ in $\R^N$.

If $\omega \in A_2(\R^N)$, we call it an $A_2$-weight, $C_{2,\omega}$ in \eqref{A_pclass} is the $A_2$-constant of $\omega$. 
\end{definition}

To study problem \eqref{alpha_harm_intro} we define the weighted Sobolev space
\begin{equation}
  \label{HL10}
  \HL(y^{\alpha},\C) = \left\{ w \in H^1(y^\alpha,\C): w = 0 \textrm{ on } \partial_L \C\right\}.
\end{equation}
As \cite[(2.21)]{NOS} shows, the following \emph{weighted Poincar\'e inequality} holds:
\begin{equation}
\label{Poincare_ineq}
\| w \|_{L^2(y^{\alpha},\C)} \lesssim \| \nabla w \|_{L^2(y^{\alpha},\C)},
\quad \forall w \in \HL(y^{\alpha},\C),
\end{equation}
then the seminorm on $\HL(y^{\alpha},\C)$ is equivalent to \eqref{wH1norm}. For $w \in H^1(y^{\alpha},\C)$, $\tr w$ denotes its trace onto $\Omega \times \{ 0 \}$ which satisfies
\cite[Proposition 2.5]{NOS}
\begin{equation}
\label{Trace_estimate}
\tr \HL(y^\alpha,\C) = \Hs,
\qquad
  \|\tr w\|_{\Hs} \leq C_{\tr} \| w \|_{\HLn(y^\alpha,\C)}.
\end{equation}

Let us now describe the Caffarelli-Silvestre extension result
for bounded domains; \cite{CS:07,ST:10}. Given $f \in \Hsd$,
let $u \in \Hs$ solve \eqref{fl=f_bdddom}. If $\ue \in
\HL(y^{\alpha},\C)$ solves \eqref{alpha_harm_intro}, then $u = \ue(\cdot,0)$ and
\eqref{eq:identity} holds.
\subsection{A priori error analysis}
\label{sec:apriori}

We now review the main results of \cite{NOS} about the a priori error analysis of discretizations of problem \eqref{fl=f_bdddom}. 
This will also serve to make clear the limitations of this theory, thereby justifying the quest for an a posteriori error analysis. In this section we assume that
\begin{equation}
\label{reg_Omega}
 \| w \|_{H^2(\Omega)} \lesssim \| \Delta_{x'} w \|_{L^2(\Omega)}, \quad \forall w \in H^2(\Omega) \cap H^1_0(\Omega). 
\end{equation}
This holds if, for instance, the domain $\Omega$ is convex \cite{Grisvard}.

Since $\C$ is unbounded, problem \eqref{alpha_harm_intro} cannot be
approximated with finite-element-like techniques. However, since the solution of problem \eqref{alpha_harm_intro} decays exponentially in $y$  \cite[Proposition 3.1]{NOS}, by truncating $\C$ to $\C_\Y$ and setting a homogeneous Dirichlet condition on $y = \Y$, we only incur in an exponentially small error in terms of $\Y$ \cite[Theorem 3.5]{NOS}.
If
\[
  \HL(y^{\alpha},\C_\Y) := \left\{ v \in H^1(y^\alpha,\C_\Y): v = 0 \text{ on }
    \partial_L \C_\Y \cup \Omega \times \{ \Y\} \right\},
\]
then the aforementioned problem reads: find $v \in \HL(y^{\alpha}, \C_\Y)$ such that
\begin{equation}
\label{alpha_harmonic_extension_weak_T}
  \int_{\C_\Y} y^\alpha \nabla v \nabla \phi
  = d_s \langle f, \tr \phi\rangle, \quad \forall v \in \HL(y^{\alpha},\C_\Y).
\end{equation}
Here $\langle \cdot, \cdot \rangle$ is the duality pairing between $\Hsd$ and $\Hs$; see \eqref{Trace_estimate}.

If $\ue$ and $v$ denote the solution of \eqref{alpha_harm_intro} and \eqref{alpha_harmonic_extension_weak_T}, respectively, then \cite[Theorem 3.5]{NOS} provides the following exponential estimate
\begin{equation*}
  \| \nabla(\ue - v) \|_{L^2(y^{\alpha},\C )} \lesssim e^{-\sqrt{\lambda_1} \Y/4} \| f\|_{\Hsd}.
\end{equation*}

To study the finite element discretization of \eqref{alpha_harmonic_extension_weak_T} we must understand the regularity of $\ue$ and $v$. As \cite[Theorem 2.7]{NOS} reveals, the second order regularity of $\ue$ is much worse in the extended direction, namely
\begin{align}
    \label{reginx}
  \| \Delta_{x'} \ue\|_{L^2(y^{\alpha},\C)} + 
  \| \partial_y \nabla_{x'} \ue \|_{L^2(y^{\alpha},\C)}
  & \lesssim \| f \|_{\Ws}, \\
\label{reginy}
  \| \ue_{yy} \|_{L^2(y^{\beta},\C)} &\lesssim \| f \|_{L^2(\Omega)},
\end{align}
with $\beta > 2\alpha + 1$. Thus, \emph{graded} meshes in the extended variable $y$ play a fundamental role. Estimates \eqref{reginx}--\eqref{reginy} motivate the construction of a mesh over $\C_{\Y}$ as follows. We consider a partition $\mathcal{I}_\Y$ of the interval $[0,\Y]$ with points
\begin{equation}
\label{graded_mesh}
  y_k = k^\gamma M^{-\gamma} \Y, \quad k=0,\dots,M,
\end{equation}
where $\gamma > 3/(1-\alpha)=3/(2s)$. Let $\T_\Omega = \{K\}$ be a conforming and shape regular mesh of $\Omega$, where $K \subset \R^n$ is an element that is isoparametrically equivalent either to the unit cube $[0,1]^n$ or the unit simplex in $\R^n$. The collection of these triangulations is denoted by $\Tr_\Omega$.  We construct the mesh $\T_{\Y}$ of $\C_\Y$ as the tensor product of $\T_\Omega$ and $\mathcal{I}_\Y$.
By construction $\T_\Y$ is anisotropic in the extended variable (\cf \cite{DL:05,NOS}). The set of all triangulations of $\C_\Y$ obtained with this procedure is $\Tr$. 

For $\T_{\Y} \in \Tr$, we define
the finite element space 
\begin{equation}
\label{eq:FESpace}
  \V(\T_\Y) = \left\{
            W \in C^0( \bar{\C}_\Y): W|_T \in \mathcal{P}_1(K) \otimes \mathbb{P}_1(I) \ \forall T \in \T_\Y, \
            W|_{\Gamma_D} = 0
          \right\},
\end{equation}
where $\Gamma_D = \partial_L \C_{\Y} \cup \Omega \times \{ \Y\}$ is
the Dirichlet boundary. The set $\mathcal{P}_1(K)$ is either the space
$\mathbb{P}_1(K)$ of polynomials of total degree at most $1$, when the
base $K$ of an element $T = K \times I$ is simplicial, or the space of polynomials
$\mathbb{Q}_1(K)$ of degree not larger than $1$ in each variable provided $K$ is an
$n$-rectangle. We also define $\U(\T_{\Omega})=\tr \V(\T_{\Y})$.
Note that $\#\T_{\Y} = M \, \# \T_\Omega$, and that $\# \T_\Omega \approx M^n$ implies $\#\T_\Y \approx M^{n+1}$.

The Galerkin approximation of \eqref{alpha_harmonic_extension_weak_T} is the function $V_{\T_{\Y}} \in \V(\T_{\Y})$ such that
\begin{equation}
\label{harmonic_extension_weak}
  \int_{\C_\Y} y^{\alpha}\nabla V_{\T_{\Y}} \nabla W = 
d_s \langle f, \textrm{tr}_{\Omega} W \rangle,
  \quad \forall W \in \V(\T_{\Y}).
\end{equation}
Existence and uniqueness of $V_{\T_{\Y}}$ immediately follows from
$\V(\T_\Y) \subset \HL(y^{\alpha},\C_\Y)$ and the Lax-Milgram
Lemma. It is trivial also to obtain a best approximation result
\emph{\`a la} Cea. This reduces the numerical analysis of
\eqref{harmonic_extension_weak} to a question in approximation theory
which in turn can be answered by the study of piecewise polynomial
interpolation in Muckenhoupt weighted Sobolev spaces; see
\cite{NOS,NOS2}. Exploiting the Cartesian structure of the mesh we can
handle the anisotropy, construct a quasi interpolant $\Pi_{\T_\Y} :
L^1(\C_\Y) \to \V(\T_\Y)$, thus extending \cite{DL:05}, and obtain
\begin{align*}
  \| v - \Pi_{\T_\Y} v \|_{L^2(y^\alpha,T)} & \lesssim 
    h_K  \| \nabla_{x'} v\|_{L^2(y^\alpha,S_T)} + h_I \| \partial_y v\|_{L^2(y^\alpha,S_T)}, \\
  \| \partial_{x_j}(v - \Pi_{\T_\Y} v) \|_{L^2(y^\alpha,T)} &\lesssim
    h_K  \| \nabla_{x'} \partial_{x_j} v\|_{L^2(y^\alpha,S_T)} + h_I \| \partial_y \partial_{x_j} v\|_{L^2(y^\alpha,S_T)},
\end{align*}
with $j=1,\ldots,n+1$; see \cite[Theorems 4.6--4.8]{NOS} and \cite{NOS2}. However, since \eqref{reginy} implies $\ue \notin H^2(y^{\alpha},\C_{\Y})$ the second estimate is not meaningful for $j=n+1$. We must measure the regularity of $\ue_{yy}$ with a stronger weight and compensate with a graded mesh. This makes anisotropic estimates essential. These considerations allow us to obtain (\cite[Theorem 5.4]{NOS} and \cite[Corollary 7.11]{NOS}):

\begin{theorem}[a priori error estimate]
\label{TH:fl_error_estimates}
Let $\T_\Y \in \Tr$ and $\V(\T_\Y)$ be defined by \eqref{eq:FESpace}. If $V_{\T_\Y} \in \V(\T_\Y)$ solves \eqref{harmonic_extension_weak}, then we have
\begin{equation*}
\label{optimal_rate}
  \| \ue - V_{\T_\Y} \|_{\HLn(y^\alpha,\C)} \lesssim
|\log(\# \T_{\Y})|^s(\# \T_{\Y})^{-1/(n+1)} \|f \|_{\mathbb{H}^{1-s}(\Omega)},
\end{equation*}
where $\Y \approx \log(\# \T_{\Y})$. Alternatively, if $u$ denotes the solution of \eqref{fl=f_bdddom}, then
\begin{equation*}
\| u - \tr V_{\T_\Y} \|_{\Hs} \lesssim
|\log(\# \T_{\Y})|^s(\# \T_{\Y})^{-1/(n+1)} \|f \|_{\mathbb{H}^{1-s}(\Omega)}.
\end{equation*}
\end{theorem}

\begin{remark}[domain and data regularity]
\label{rm:dom_and_data2}
The results of Theorem~\ref{TH:fl_error_estimates} hold only if $f \in \mathbb{H}^{1-s}(\Omega)$ 
and the domain $\Omega$ is such that \eqref{reg_Omega} holds.
\end{remark}

\begin{remark}[quasi-uniform meshes]
\label{thm:quasiuniform}
Let $\ue$ solve \eqref{alpha_harm_intro}, and $V_{\T_{\Y}}$ be the solution of \eqref{harmonic_extension_weak}, constructed over a quasi-uniform mesh of size $h_{\T}$. If $f \in \Ws$ and $\Y \approx |\log h_{\T} |$, then for all $\varepsilon >0$
\begin{equation}
  \label{sub_optimal}
   \| \nabla(\ue - V_{\T_{\Y}} ) \|_{L^2(y^\alpha, \C_{\Y})}  \lesssim { h_{\T}^{s-\varepsilon} }\| f \|_{\Ws},
\end{equation}
where the hidden constant blows up if $\varepsilon$ tends to $0$.
\end{remark}

\begin{remark}[Case $s = \srn$]
If $s=\srn$, then there is no weight in \eqref{sub_optimal} and we obtain the optimal estimate
\[
   \| \nabla(\ue - V_{\T_{\Y}} ) \|_{L^2(\C_{\Y})}  
   \lesssim h_{\T}\| f \|_{H^{1/2}_{00}(\Omega)}.
\]
\end{remark}

\subsection{Numerical illustrations}
\label{sec:numerics}

Let us present several numerical examples. These have been implemented with the \texttt{deal.II} library \cite{dealii,dealii2}. Integrals are evaluated with Gaussian quadratures of sufficiently high order and linear systems are solved using CG with ILU preconditioner and the exit criterion being that the $\ell^2$-norm of the residual is less than $10^{-12}$. 

\subsubsection{Quasi-uniform meshes}

\begin{figure}[h!]
  \centering
  \includegraphics[width=0.45\textwidth]{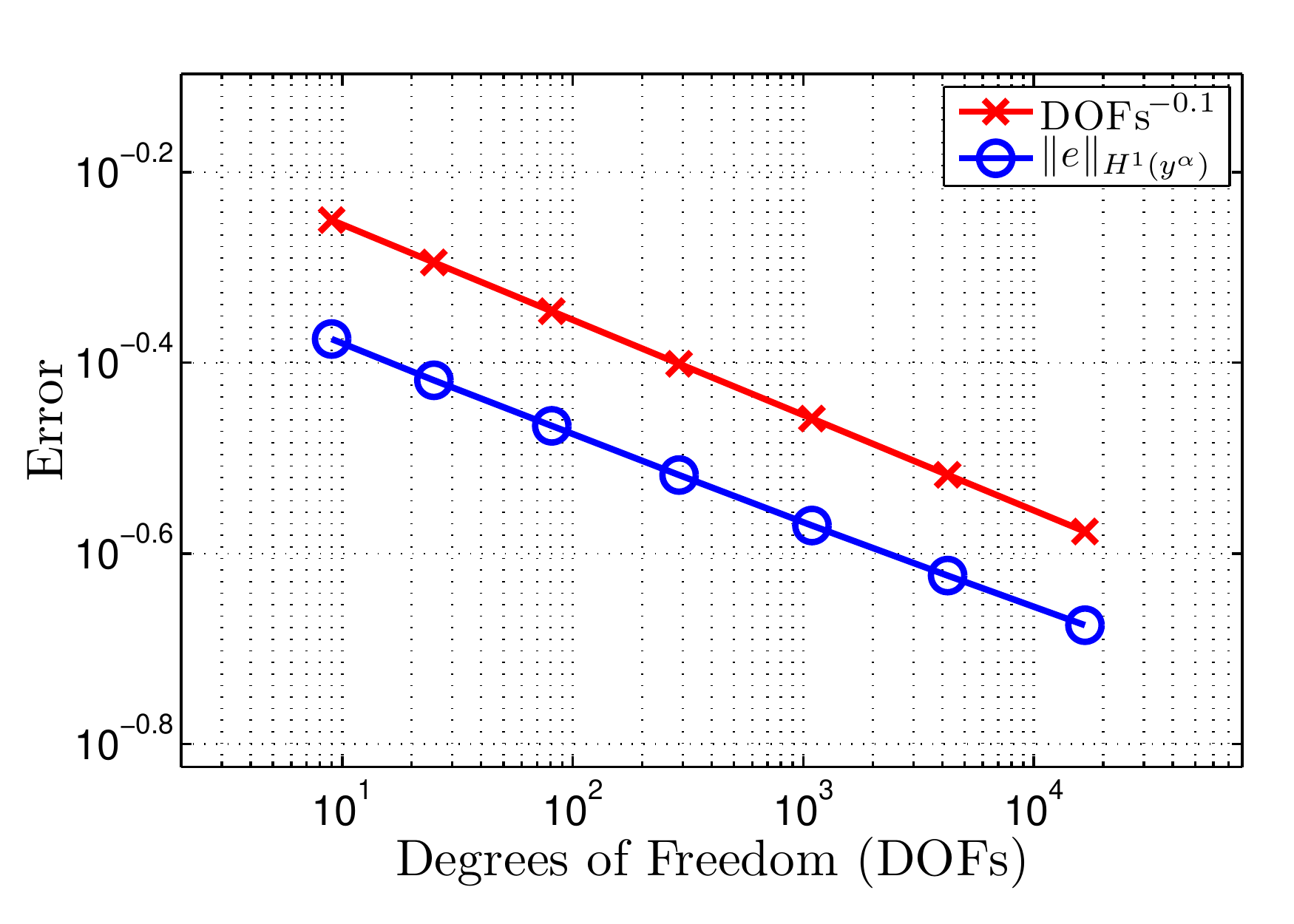}
  \vskip-0.3cm
  \caption{Rate of convergence $\# (\T_{\Y})^{-s/(n+1)}$ for quasi-uniform meshes: $s=0.2$, $n=1$.}
  \label{fig:0.2_uniform}
\end{figure}

According to \cite{NOS}, $\partial_y \ue \approx y^{-\alpha}$, which
formally implies $\partial_y \ue \in H^{s-\varepsilon}(y^{\alpha},\C)$
provided $f\in\Ws$. In this sense \eqref{sub_optimal} seems sharp with
respect to regularity even though it does not exhibit the optimal
rate. We verify this with a simple numerical illustration
in Figure~\ref{fig:0.2_uniform}.

\noindent
Let $n=1$, $\Omega = (0,1)$, $s=0.2$, $f = \pi^{2s} \sin(\pi x)$, then $u(x)=\sin(\pi x)$, and the solution $\ue$ to \eqref{alpha_harm_intro} is
$
   \ue(x,y) = 2^{1-s}\pi^{s}\Gamma(s)^{-1} \sin (\pi x) K_{s}(\pi y),
$
where by $K_s$ we denote the modified Bessel function of the second kind; \cite[\S 2.4]{NOS}. Figure~\ref{fig:0.2_uniform} shows the rate of convergence for the  $H^1(y^\alpha,\C_\Y)$-seminorm. Estimate \eqref{sub_optimal} predicts a rate of $h_{\T}^{0.2-\varepsilon}$.
In other words the rate of convergence, when measured in terms of degrees of freedom,  is $(\# \T_{\Y})^{-0.1-\varepsilon}$, which is what Figure~\ref{fig:0.2_uniform} displays.
%
  
\subsubsection{Graded meshes: a square domain}
\label{subsub:square}

Let $\Omega = (0,1)^2$. Then
\[
  \varphi_{m,n}(x_1,x_2) = \sin(m \pi x_1)\sin(n \pi x_2),
  \quad
  \lambda_{m,n} = \pi^2 \left( m^2 + n^2 \right),
  \qquad m,n \in \mathbb{N}.
\]
If $f(x_1,x_2) = ( 2\pi^2)^{s} \sin(\pi x_1)\sin(\pi x_2)$, then
$
  u(x_1,x_2) = \sin(\pi x_1)\sin(\pi x_2),
$
by \eqref{def:second_frac} and
$
  \ue(x_1,x_2,y) = 2^{1-s/2} \pi^s \Gamma(s)^{-1} \sin(\pi x_1)\sin(\pi x_2) y^{s}K_s(\sqrt{2}\pi y)
$
\cite[(2.24)]{NOS}.

We construct a sequence of  meshes $\{\T_{\Y_k} \}_{k\geq1}$, where $\T_\Omega$ is obtained by uniform refinement and $\calI_{\Y_k}$ is given by \eqref{graded_mesh} with parameter $\gamma > 3/(1-\alpha)$. On the basis of Theorem~\ref{TH:fl_error_estimates}, the truncation parameter $\Y_k$ is chosen by
\[
  \Y_k \geq \frac{2}{\sqrt{\lambda_1}} ( \log C - \log (\# \T_{\Y_{k-1}})^{-1/3} ).
\] 
With this type of meshes,
\[
  \| u - \tr V_{\T_{\Y,k}} \|_{\Hs} \lesssim  \| \ue -V_{\T_{\Y_k}} \|_{\HLn(y^\alpha,\C)} \lesssim  {|\log( \# \T_{\Y_k} )|^s} \cdot (\# \T_{\Y_k})^{-1/3},
\]
which is near-optimal in $\ue$ but suboptimal in $u$, since we should expect (see \cite{BrennerScott})
\[
  \| u - \tr V_{\T_{\Y,k}}  \|_{\Hs} \lesssim h_{\T_\Omega}^{2-s} \lesssim (\# \T_{\Y_k} )^{-(2-s)/3}.
\] 

\begin{figure}[h!]
\label{fig:02}
\begin{center}
  \includegraphics[scale=0.3]{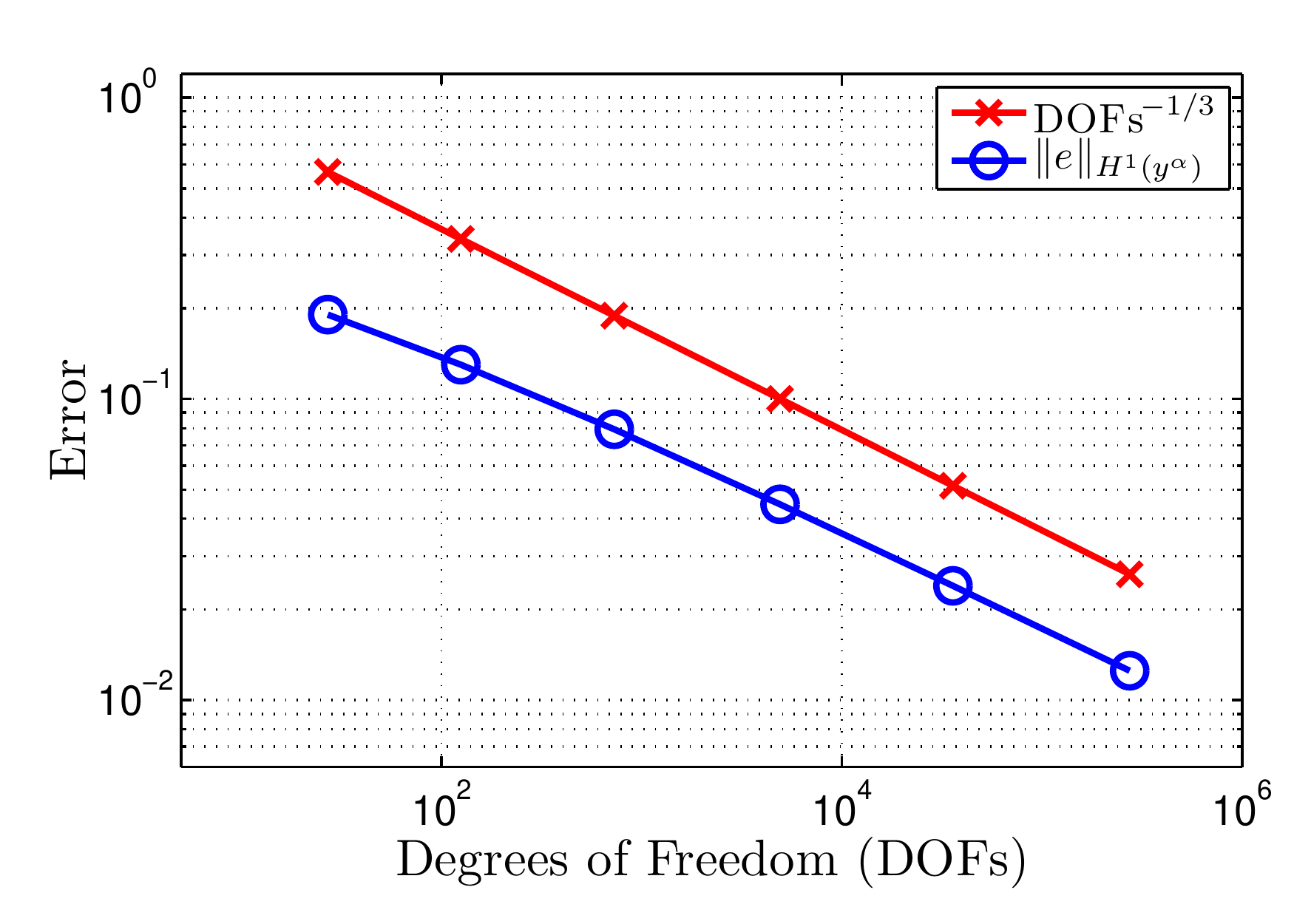}
  \includegraphics[scale=0.3]{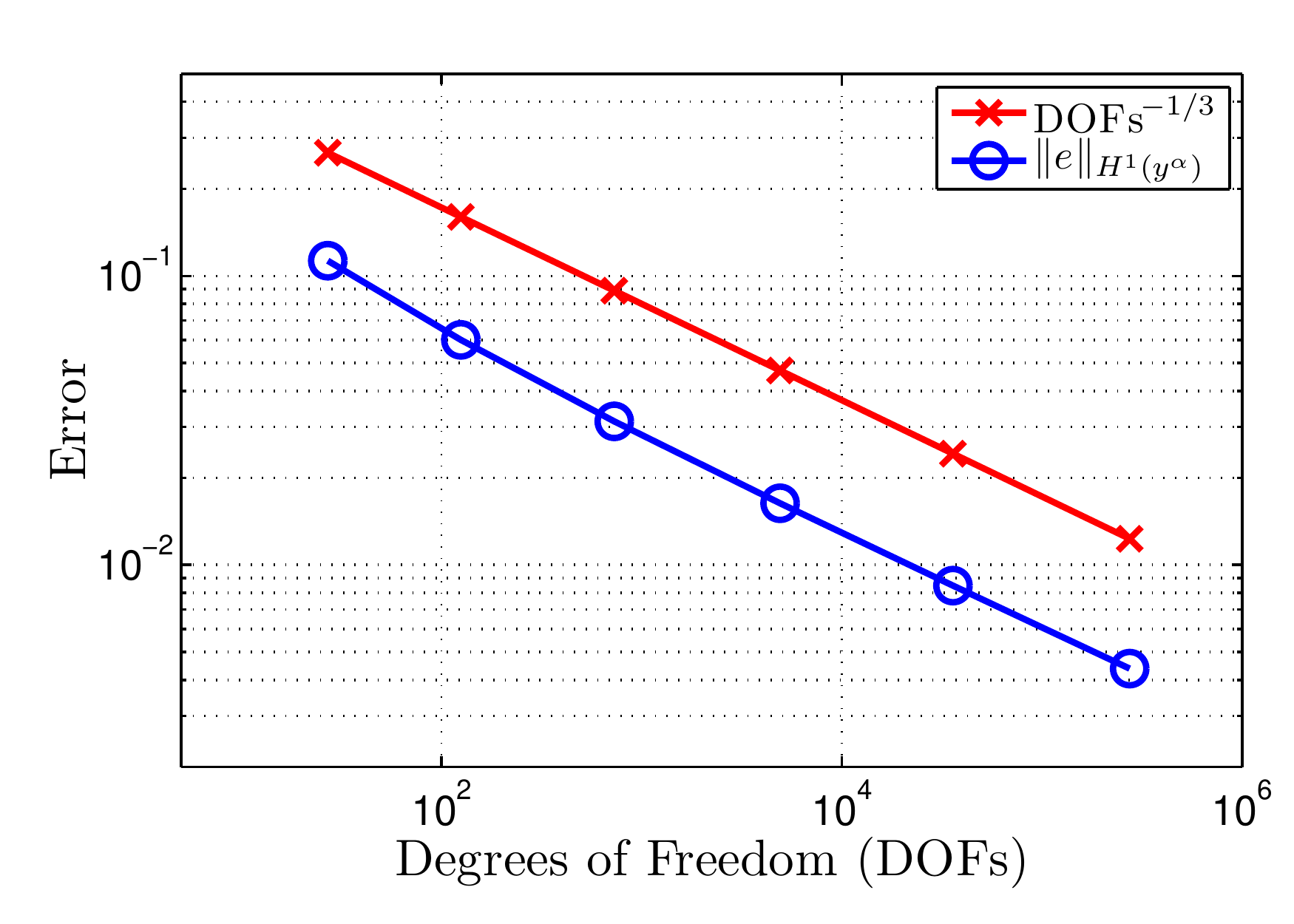}
\end{center}
\vskip-0.3cm
\caption{Computational rate of convergence for a square and graded meshes. The left panel shows the rate for $s=0.2$ and the right one 
for $s=0.8$. The experimental rate of convergence is $(\# \T_{\Y_k}
)^{-1/3}$, in agreement with Theorem~\ref{TH:fl_error_estimates}.}
\end{figure}

Figure~\ref{fig:02} shows the rates of convergence for $s=0.2$ and $s=0.8$ respectively.
In both cases, we obtain the rate given by Theorem~\ref{TH:fl_error_estimates}.

\subsubsection{Graded meshes: a circular domain}
If $\Omega = \{ |x'| \in \R^2 : |x'|<1 \}$, then
\begin{equation}
\label{circularphi}
 \varphi_{m,n}(r,\theta) = 
 J_{m}(j_{m,n}r)
\left( A_{m,n} \cos (m \theta) + B_{m,n} \sin(m \theta)\right),
\
\lambda_{m,n}=j_{m,n}^2,
\end{equation}
where $J_m$ is the $m$-th Bessel function of the first kind;
$j_{m,n}$ is the $n$-th zero of $J_m$ and
$A_{m,n}$, $B_{m,n}$ are normalization constants to ensure
$\| \varphi_{m,n} \|_{L^2(\Omega)}=1$.

If $f = ( \lambda_{1,1})^{s} \varphi_{1,1}$, then \eqref{def:second_frac} and \cite[(2.24)]{NOS} show that $u = \varphi_{1,1}$ and
\[
  \ue(r,\theta,y) = 2^{1-s}\Gamma(s)^{-1}(\lambda_{1,1})^{s/2}
         \varphi_{1,1}(r,\theta) y^{s}K_s(\sqrt{2}\pi y).
\]

\begin{figure}
\begin{center}
  \includegraphics[scale=0.3]{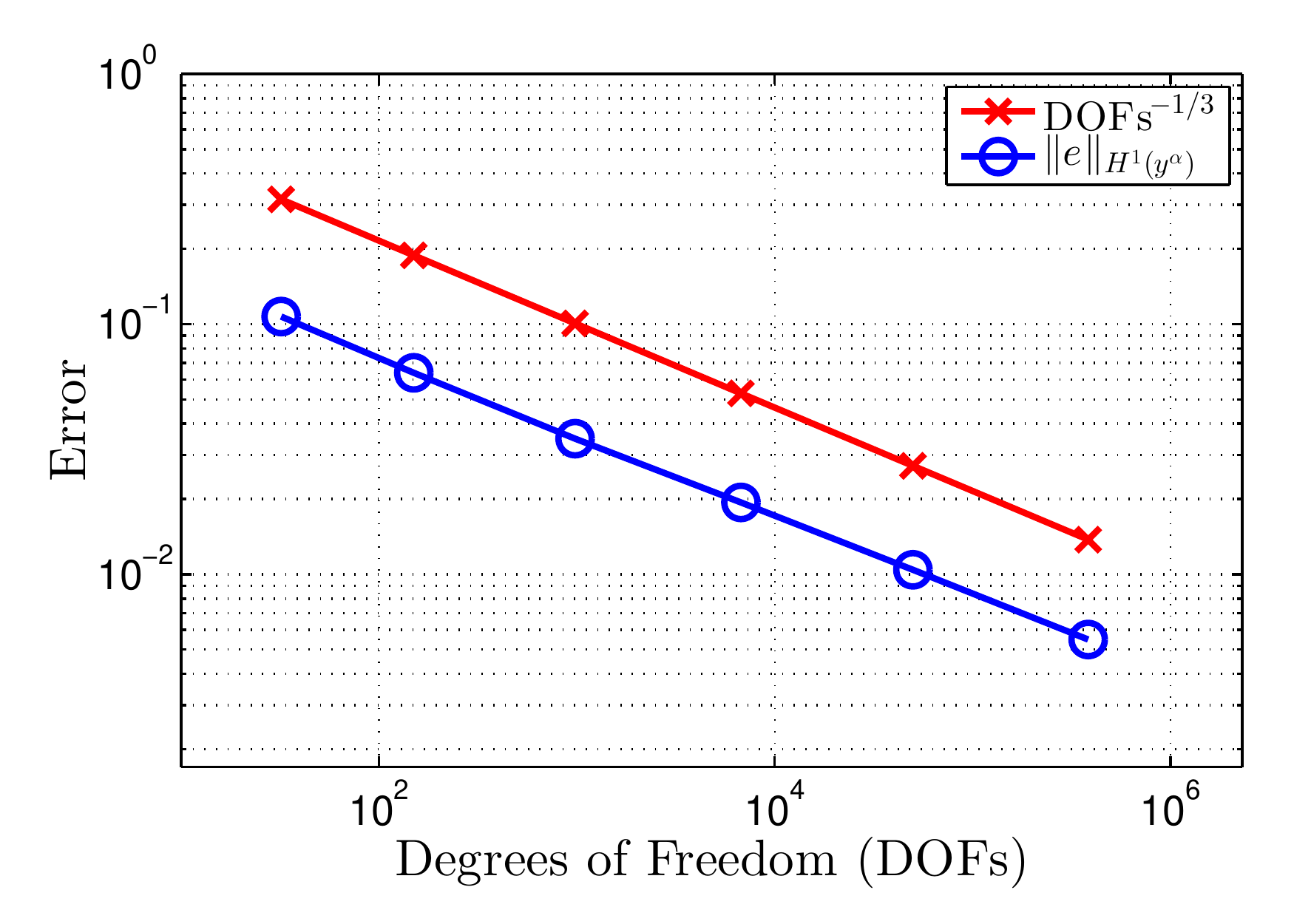}
  \includegraphics[scale=0.3]{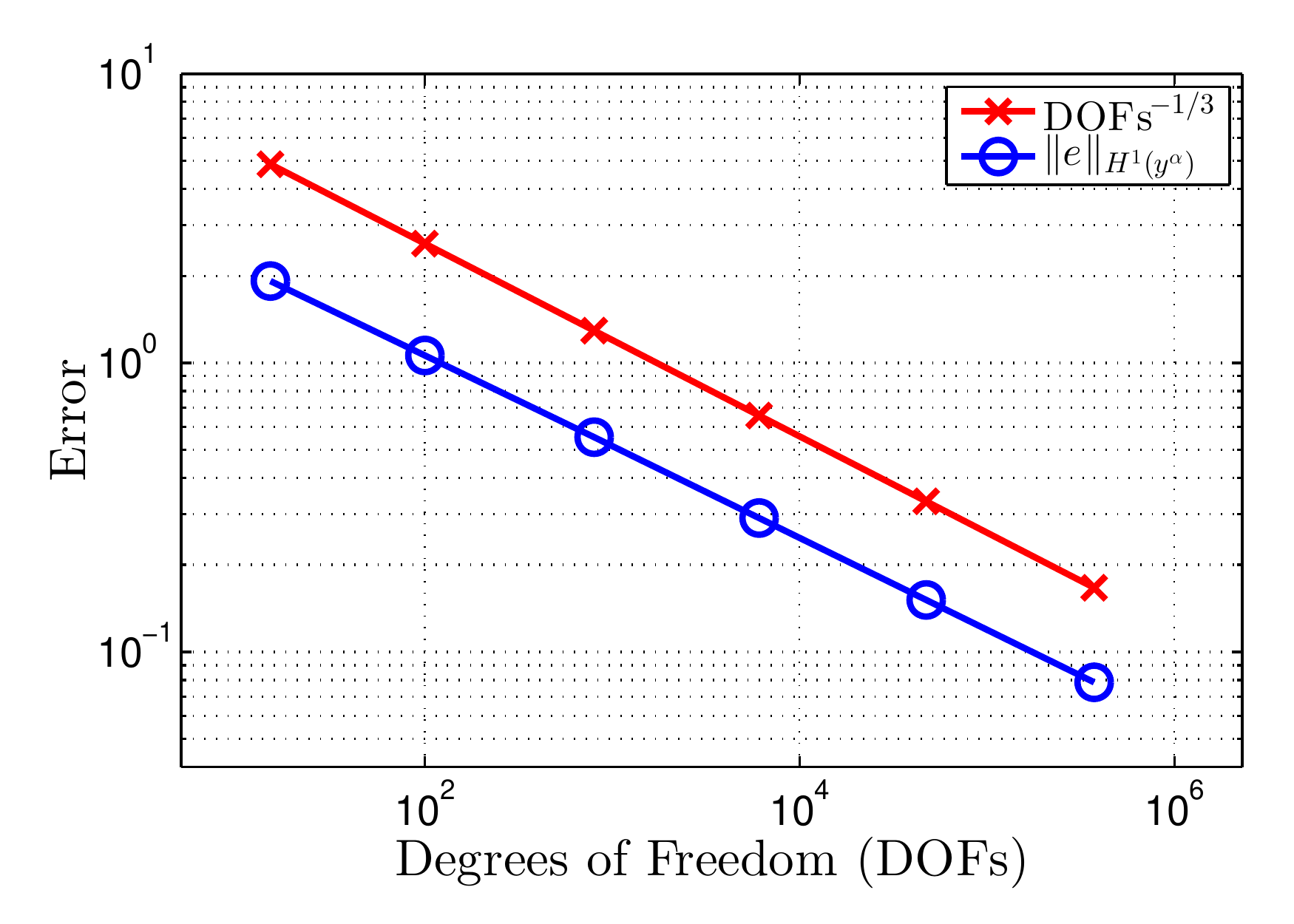}
\end{center}
\vskip-0.3cm
\caption{Computational rate of convergence for a circle with graded meshes. The left panel shows the rate for $s=0.3$ and the right one 
for $s=0.7$. The experimental rate of convergence is $(\#\T_{\Y_k} )^{-1/3}$,
in agreement with Theorem~\ref{TH:fl_error_estimates}.}
\label{fig:03}
\end{figure}

We construct a sequence of  meshes $\{\T_{\Y_k} \}_{k\geq1}$ as in \S\ref{subsub:square}. With these meshes
\begin{equation}
\label{numerical_experiment_2_ve}
  \| \ue -V_{\T_{\Y_k}} \|_{\HLn(\C,y^{\alpha})} \lesssim 
  |\log(\# \T_{\Y_k})|^s (\# \T_{\Y_k})^{-1/3},
\end{equation}
which is near-optimal. Figure~\ref{fig:03} shows the errors of $\|\ue - V_{\T_{\Y_k}} \|_{H^1(y^{\alpha},\C_{\Y_{k}})}$
for $s = 0.3$ and $s = 0.7$. The results, again, are in agreement with
Theorem~\ref{TH:fl_error_estimates}.

\section{A posteriori error analysis and adaptivity}
\label{sec:a_posteriori}
Starting with the pioneering works \cite{BabuskaMiller,BR78} in the late 1970's, adaptive finite element methods (AFEM) have been the subject of intense research. This is because AFEM yields optimal performance in situations where classical FEM cannot. In our problem we are incorporating one extra dimension and so AFEM is essential to retain computational efficiency. In addition, the a priori theory requires $f \in \Ws$ and \eqref{reg_Omega}. If either of these does not hold $\ue$ may have singularities in the $x'$-variables and exhibit fractional regularity. This would not allow us to attain the almost optimal rate of convergence of Theorem~\ref{TH:fl_error_estimates}: a quasi-uniform refinement of $\Omega$ would not result in an efficient solution technique. An adaptive loop driven by an a posteriori error estimator is \emph{fundamental} to recover optimal rates of convergence.
Let us explore this now and develop a new estimator.

\subsection{Residual estimators}

A residual error estimator uses the strong form of the local residual to measure the error. Let $T \in \T_{\Y}$ and $\nu$ be the unit outer normal to $T$. Integration by parts yields
\[
 \int_{T} y^{\alpha} \nabla V_{\T_\Y} \cdot \nabla W = \int_{\partial T} W y^{\alpha} \nabla V_{\T_\Y} \cdot \nu 
- \int_{T} \DIV(y^{\alpha} \nabla V_{\T_\Y}) W.
\]
Since $\alpha \in (-1,1)$ the boundary integral is meaningless for $y = 0$. Even the very first step in the derivation of a residual a posteriori error estimator fails! There is nothing left to do but to consider a different type of estimator.

\subsection{Local problems on stars over isotropic refinements}
\label{subsec:isotropic}

Following \cite{BabuskaMiller,MNS02} we can construct, over shape regular meshes, an error estimator based on the solution of small problems on stars. Its construction and analysis is  similar to the developments of \S\ref{sub:cylstars}, so we skip details. Under some assumptions this estimator is equivalent to the error up to data oscillation; see \cite{NOS3}. We designed an adaptive algorithm driven by this error estimator on shape regular meshes \cite{MR2875241,MNS02}, and we illustrate its performance with a simple but revealing numerical example.

\begin{figure}[h!]
  \begin{center}
    \includegraphics[width=0.43\textwidth]{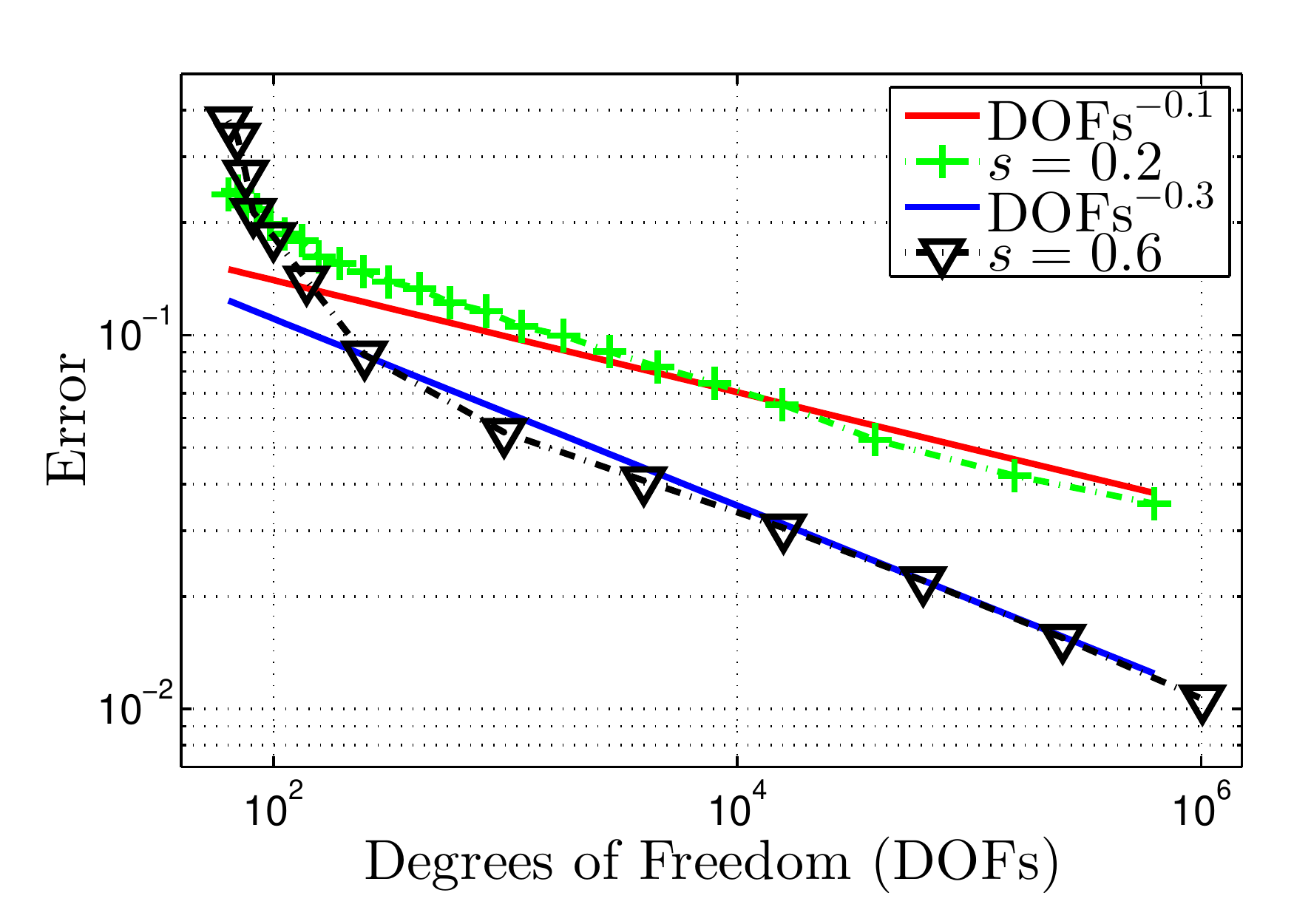}
  \end{center}
  \vskip-0.3cm
  \caption{Computational rate of convergence $\# (\T_{\Y})^{-s/2}$ for an isotropic adaptive algorithm
   for $n=1$, $s=0.2$ and $s=0.6$. }
\label{fig:s0.2isotropic}
\end{figure}

Let $\Omega = (0,1)$ and $s \in (0,1)$. If $f(x') = \pi^{2s} \sin(\pi x')$, then $u(x')=\sin(\pi x')$, and
$
   \ue(x',y) = 2^{1-s}\pi^{s}\Gamma(s)^{-1} \sin (\pi x') K_{s}(\pi y)
$ solves \eqref{alpha_harm_intro}.
We point out that we are solving a two dimensional problem so the optimal rate we expect is $\mathcal{O}( \# \, \T_\Y^{-0.5})$. Figure~\ref{fig:s0.2isotropic} shows the experimental rate of convergence for $s = 0.2$ and $s = 0.6$ which, as we see, is
$
 \mathcal{O}(\# \T_\Y^{-s/2 })
$
and coincides with the suboptimal one obtained with quasi-uniform
refinement \cite[\S5.1]{NOS}. These numerical experiments show that
adaptive isotropic refinement cannot be optimal, thus justifying the
need to introduce \emph{cylindrical stars} together with a new
anisotropic error estimator, which will treat the $x'$-coordinates and
the extended direction $y$ separately.


\subsection{Local problems over cylindrical stars}
\label{sub:cylstars}
Inspired by \cite{BabuskaMiller,MNS02}, we deal with the anisotropy of the mesh in the extended variable $y$ and the coefficient $y^{\alpha}$ by considering local problems on cylindrical stars. The solutions of these local problems allow us to define an anisotropic a posteriori error estimator which, under certain assumptions, is equivalent to the error up to data oscillation terms.

Given a node $z$ on the mesh $\T_{\Y}$, we exploit the tensor product structure of $\T_{\Y}$, and  we write $z = (z',z'')$ where $z'$ and $z''$ are nodes on the meshes 
$\T_{\Omega}$ and $\mathcal{I}_{\Y}$ respectively. For $K \in \T_{\Omega}$, we denote by $\N(K)$ and $\Nin(K)$ the set of nodes and interior nodes of $K$, respectively. We set 
\[
\N(\T_{\Omega}) = \bigcup_{K \in \T_\Omega} \N(K), \qquad \Nin(\T_{\Omega}) = \bigcup_{K \in \T_\Omega} \Nin(K).
\]
The \emph{star} around $z'$ is
$
  S_{z'} = \bigcup_{K \ni z'} K \subset \Omega,
$
and the \emph{cylindrical star} around $z'$ is
\[
  \C_{z'} := \bigcup\left\{ T \in \T_\Y : T = K \times I,\ K \ni z'  \right\}= S_{z'} \times (0,\Y) \subset \C_{\Y}.
\]

For each node $z' \in \N(\T_\Omega)$ we define the discrete space
\begin{align*}
  \mathbb{W}(\C_{z'}) &= 
    \left\{
      W \in C( \bar{\C}_\Y): W|_T \in \mathcal{P}_2(K) \otimes \mathbb{P}_2(I) \ \forall T = K \times I \in \C_{z'}, 
      \right. \\
    &\left.
      W|_{\partial \C_{z'} \setminus \Omega \times \{ 0\} } = 0
    \right\},
\end{align*}
where, if $K$ is a quadrilateral, $\mathcal{P}_2(K)=\mathbb{Q}_2(K)$ --- the space of polynomials of degree not larger than $2$ in each variable. If $K$ is a simplex, $\mathcal{P}_2(K)=\mathbb{P}_2(K) \oplus \mathbb{B}(K)$ where $\mathbb{P}_2(K)$ is the space of polynomials of total degree at most $2$, and $\mathbb{B}(K)$ is the space spanned by a local cubic bubble function.

For each cylindrical star $\C_{z'}$ we define $\eta_{z'} \in \mathbb{W}(\C_{z'})$ to be the solution of
\begin{equation*}
  \int_{\C_{z'}} y^{\alpha} \nabla\eta_{z'} \nabla W =
  d_s \langle f, \tr W \rangle
   -  \int_{\C_{z'}} y^{\alpha} \nabla V_{\T_{\Y}} \nabla W,
   \quad \forall W \in \mathbb{W}(\C_{z'}).
\end{equation*}
We finally define the local indicators $\E_{z'}$ and global error estimators
$\E_{\T_{\Y}}$ as follows:
\begin{equation}
\label{comp_global_estimator}
  \E_{z'} = \interleave \eta_{z'} \! \!  \interleave_{\C_{z'}} \qquad
  \E_{\T_{\Y}}^2 = \sum_{z' \in \N(\T_{\Omega}) } \E_{z'}^2.
\end{equation}
We can obtain a local lower bound for the error without any oscillation term and free of any constant.  
  
\begin{theorem}[localized lower bound]
\label{th:dis_lower_bound}
Let $v \in \HL(y^\alpha,\C_\Y)$ and $V_{\T_\Y} \in \V(\T_\Y)$ solve 
\eqref{alpha_harmonic_extension_weak_T} and \eqref{harmonic_extension_weak} respectively. 
Then, for any $z' \in \N(\T_{\Omega})$, we have 
\begin{equation*}
\E_{z'} \leq\| \nabla(v-V_{\T_{\Y}})\|_{L^2(y^{\alpha},\C_{z'})}.
\end{equation*}
\end{theorem}

For $z' \in \N(\T_\Omega)$, we define the local \emph{data oscillation} as
\begin{equation}
\label{eq:defoflocosc}
  \osc_{z'}(f)^2 := d_s h_{z'}^{2s} \| f - f_{z'} \|_{L^2(S_{z'})}^2,
  \quad
  f_{z'}|_K := \fint_{K} f,
\end{equation}
whence the global data oscillation reads
\begin{equation*}
  \osc_{\T_\Omega}(f)^2 := \sum_{z' \in \N(\T_\Omega)} \osc_{z'}( f)^2 .
\end{equation*}
To mark elements for refinement we use the \emph{total error indicator}
\begin{equation}
\label{total_error}
\tau_{\T_{\Omega}} (V_{\T_{\Y}}, S_{z'}) := \left( \E_{z'}^2 + \osc_{z'}(f)^2 \right)^{1/2} \quad 
\forall z' \in \N(\T_{\Omega}).
\end{equation}
Let $\mathscr{K}_{\T_{\Omega}} = \{ S_{z'} : z' \in \N(\T_\Omega)\}$ and, for any $\mathscr{M} \subset \mathscr{K}_{\T_{\Omega}}$, we set 
\begin{equation}\label{total_est}
\tau_{\T_{\Omega}} (V_{\T_{\Y}}, \mathscr{M} ) := 
 \left( \sum_{S_{z'} \in \mathscr{M} } \tau_{\T_{\Omega}}(V_{\T_{\Y}}, S_{z'})^2 \right)^{1/2}.
\end{equation}

Under certain assumptions we can bound the error by the estimator, up to oscillation terms; see \cite{CNOS2} for details.

\begin{theorem}[global upper bound]
\label{th:upboundcomp}
Let $v \in \HL(\C_\Y,y^\alpha)$ and $V_{\T_\Y} \in \V(\T_\Y)$ 
solve \eqref{alpha_harmonic_extension_weak_T} and \eqref{harmonic_extension_weak}, respectively. 
The total error estimator $\tau_{\T_{\Omega}} (V_{\T_{\Y}},\mathscr{K}_{\T_{\Omega}})$, defined in 
\eqref{total_est} satisfies
\begin{equation*}
 \| \nabla (v - V_{\T_\Y}) \|_{L^2(y^\alpha,\C_\Y)} \lesssim 
\tau_{\T_{\Omega}}(V_{\T_\Y},\mathscr{K}_{\T_{\Omega}}).
\end{equation*}
\end{theorem}

\subsection{Numerical experiments}
We now illustrate the performance of the a posteriori error estimator \eqref{comp_global_estimator}. We use an almost standard adaptive loop
\begin{equation}
 \textsf{\textup{SOLVE}} \rightarrow \textsf{\textup{ESTIMATE}} \rightarrow \textsf{\textup{MARK}} \rightarrow \textsf{\textup{REFINE}}.
\label{afem}
\end{equation}
where the modules in \eqref{afem} are as follows:
\begin{enumerate}[$\bullet$]
\item \textsf{\textup{SOLVE}}: Given a mesh $\T_{\Y}$ we compute $V_{\T_\Y}$, the solution of \eqref{harmonic_extension_weak}.

\item \textsf{\textup{ESTIMATE}}: Given $V_{\T_\Y}$ we calculate 
the local error indicators \eqref{comp_global_estimator} and the local oscillations \eqref{eq:defoflocosc}
to construct  the total error indicator \eqref{total_error}:
\[
 \left\{ \tau_{\T_{\Omega}} (V_{\T_{\Y}}, S_{z'}) \right\}_{ S_{z'} \in \mathscr{K}_{\T_\Omega}} = \textsf{\textup{ESTIMATE}}(V_{\T_{\Y}}, \T_{\Y}).
\]

\item \textsf{\textup{MARK}}: Using D\"{o}rfler marking \cite{MR1393904} with parameter $0< \theta \leq 1$ we select a set 
\[
  \mathscr{M} = \textsf{\textup{MARK}}( \left\{ \tau_{\T_{\Omega}} (V_{\T_{\Y}}, S_{z'}) \right\}_{ S_{z'} \in \mathscr{K}_{\Omega}}, V_{\T_{\Y}} ) 
  \subset \mathscr{K}_{\T_{\Omega}}
\]
of minimal cardinality that satisfies
$
  \tau_{\T_{\Omega}}( V_{\T_{\Y}},\mathscr{M})  \geq \theta \tau_{\T_{\Omega}}(V_{\T_{\Y}},  \mathscr{K}_{\T_{\Omega}}). 
$

\item \textsf{\textup{REFINE}}: We generate a new mesh $\T_\Omega'$ by bisecting all the elements 
$K \in \T_{\Omega}$ contained in $\mathscr{M}$ based on newest-vertex bisection method; \cite{NV,NSV:09}. 
We choose the truncation parameter as $\Y = 1 + \tfrac{1}{3}\log(\# \T_{\Omega}')$ \cite[Remark 5.5]{NOS}.
We set $M \approx (\# \T_\Omega')^{1/n}$ and construct $\mathcal{I}_\Y'$ by the rule \eqref{graded_mesh}. The new mesh
$
  \T_\Y' = \textsf{\textup{REFINE}}(\mathscr{M})
$
is obtained as the tensor product of $\T_\Omega'$ and $\mathcal{I}_\Y'$.
\end{enumerate}

\subsubsection{Smooth but incompatible data}
\label{sub:smoothincompatible}
\begin{figure}[h!]
\centering
\includegraphics[width=0.4\textwidth]{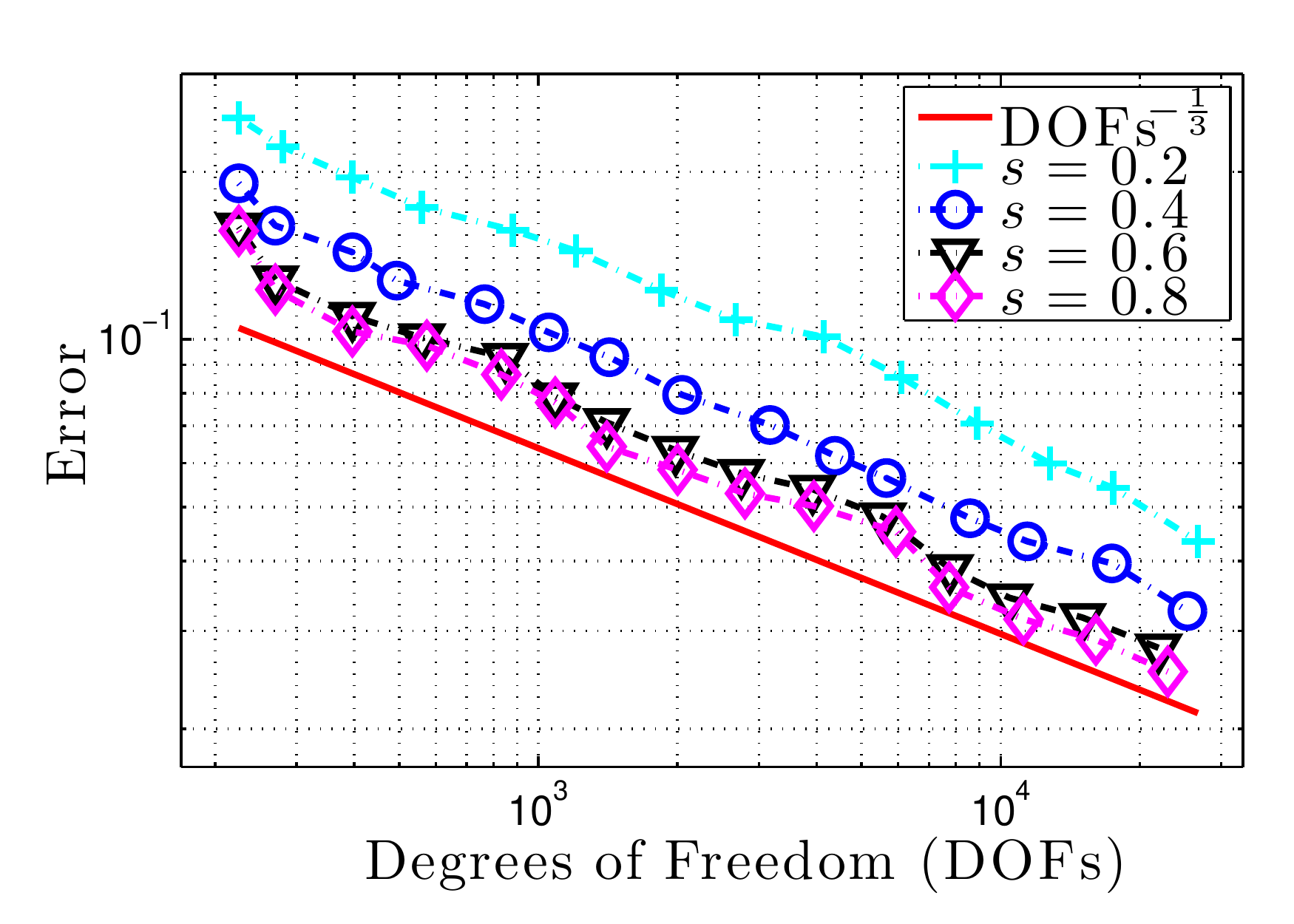}
\hfil
\includegraphics[width=0.4\textwidth]{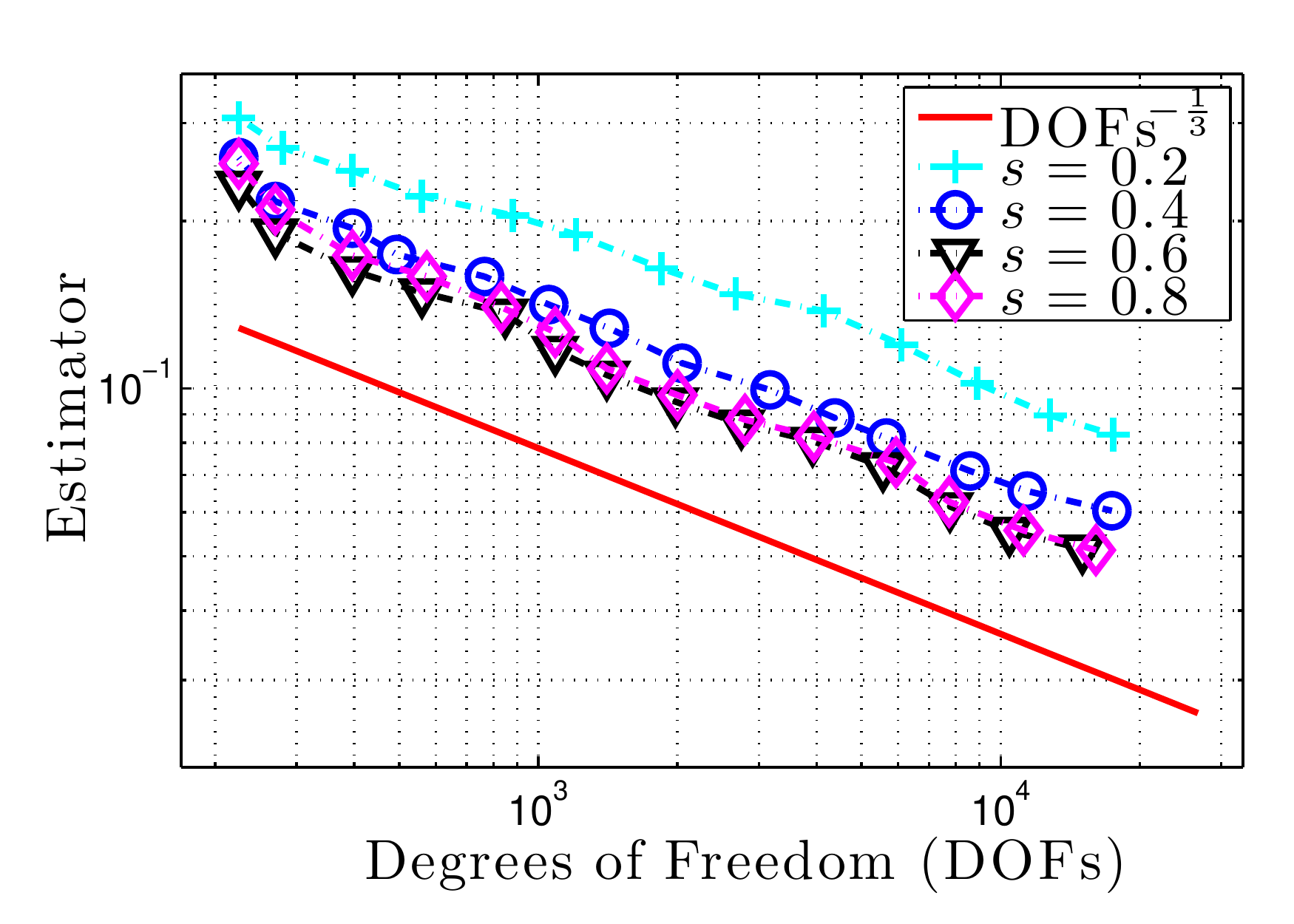}
\caption{Computational rate of convergence of AFEM for \S~\ref{sub:smoothincompatible}, $n=2$ and $s=0.2$, $0.4$, $0.6$ and $s=0.8$. The left panel shows the the error vs.~the number of degrees of freedom, the right one the total error indicator. We recover the optimal rate $\# (\T_{\Y})^{-1/3}$. For $s<\tfrac12$, the right hand side $f=1\notin \Ws$ and a quasiuniform mesh in $\Omega$ does not deliver the optimal rate of convergence \cite[\S6.3]{NOS}.}
\label{fig:incompatible}
\end{figure}

The example of \cite[\S6.3]{NOS} shows that
Theorem~\ref{TH:fl_error_estimates} is sharp: $f\in\Ws$ is necessary
to obtain an optimal rate of convergence with a quasiuniform mesh in
the $x'$-direction. A certain compatibility between the data and the
boundary condition is necessary. Moreover, \cite[\S6.3]{NOS} shows
how, in some simple cases, one can guess the singularity and a priori
design a mesh that captures it and recover the optimal rate of
convergence. This is not always possible and here we show that the
estimator \eqref{comp_global_estimator} automatically produces a
sequence of meshes that yield the optimal rate of
convergence. Consider $\Omega = (0,1)^2$ and $f=1$. From \eqref{def:Hs} we see that functions in $\Ws$, for $1-s>\tfrac12$, must have a vanishing trace. Therefore, if $s < \tfrac12$, $f\notin \Ws$  and Theorem~\ref{TH:fl_error_estimates} cannot be invoked. Nevertheless, as  Figure~\ref{fig:incompatible} shows, we recover the optimal rate.

\subsubsection{L-shaped domain with incompatible data}
\label{sub:Lshapedincompatible}

\begin{figure}[h!]
  \begin{center}
  \includegraphics[width=0.4\textwidth]{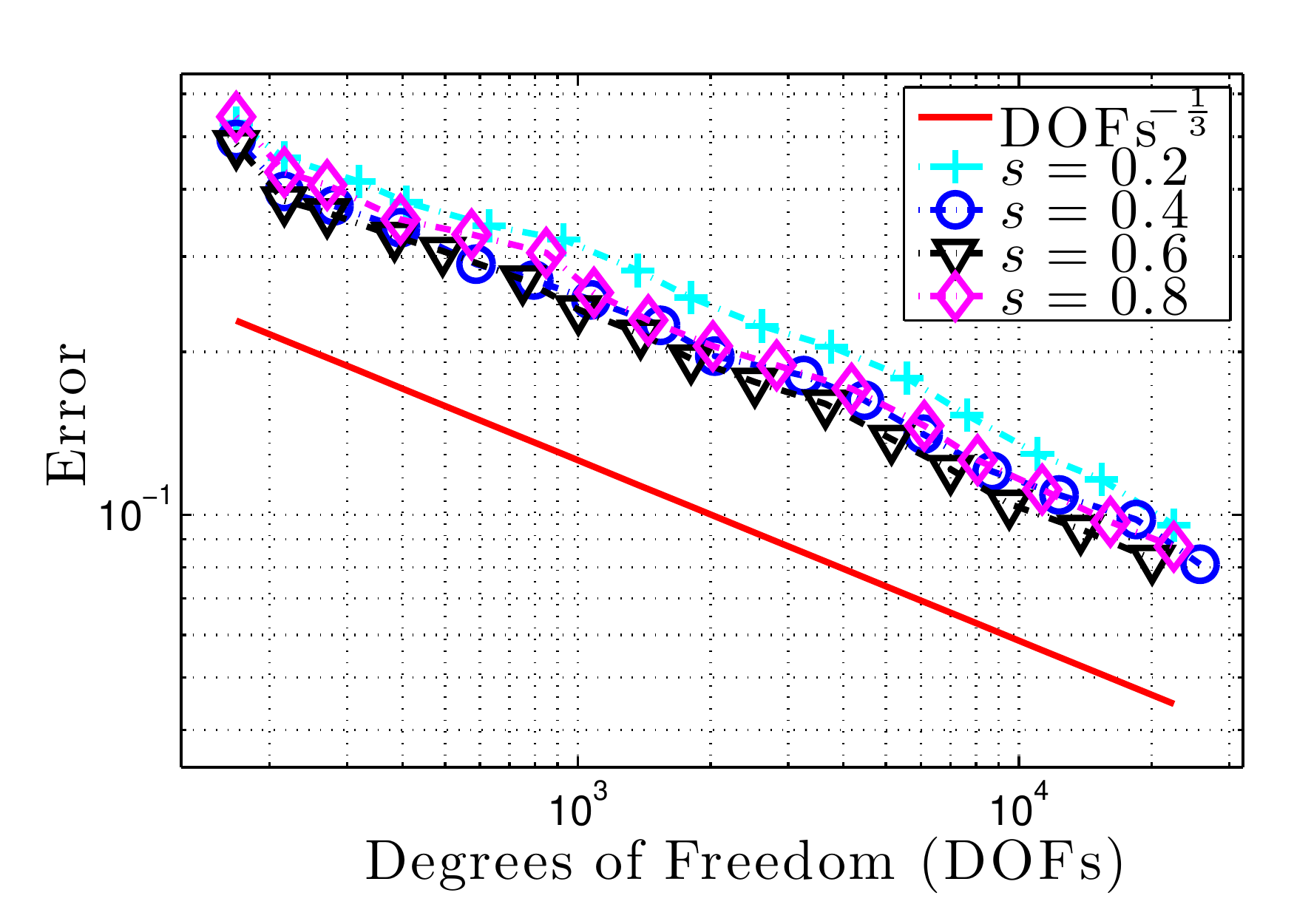}
  \hfil
  \includegraphics[width=0.4\textwidth]{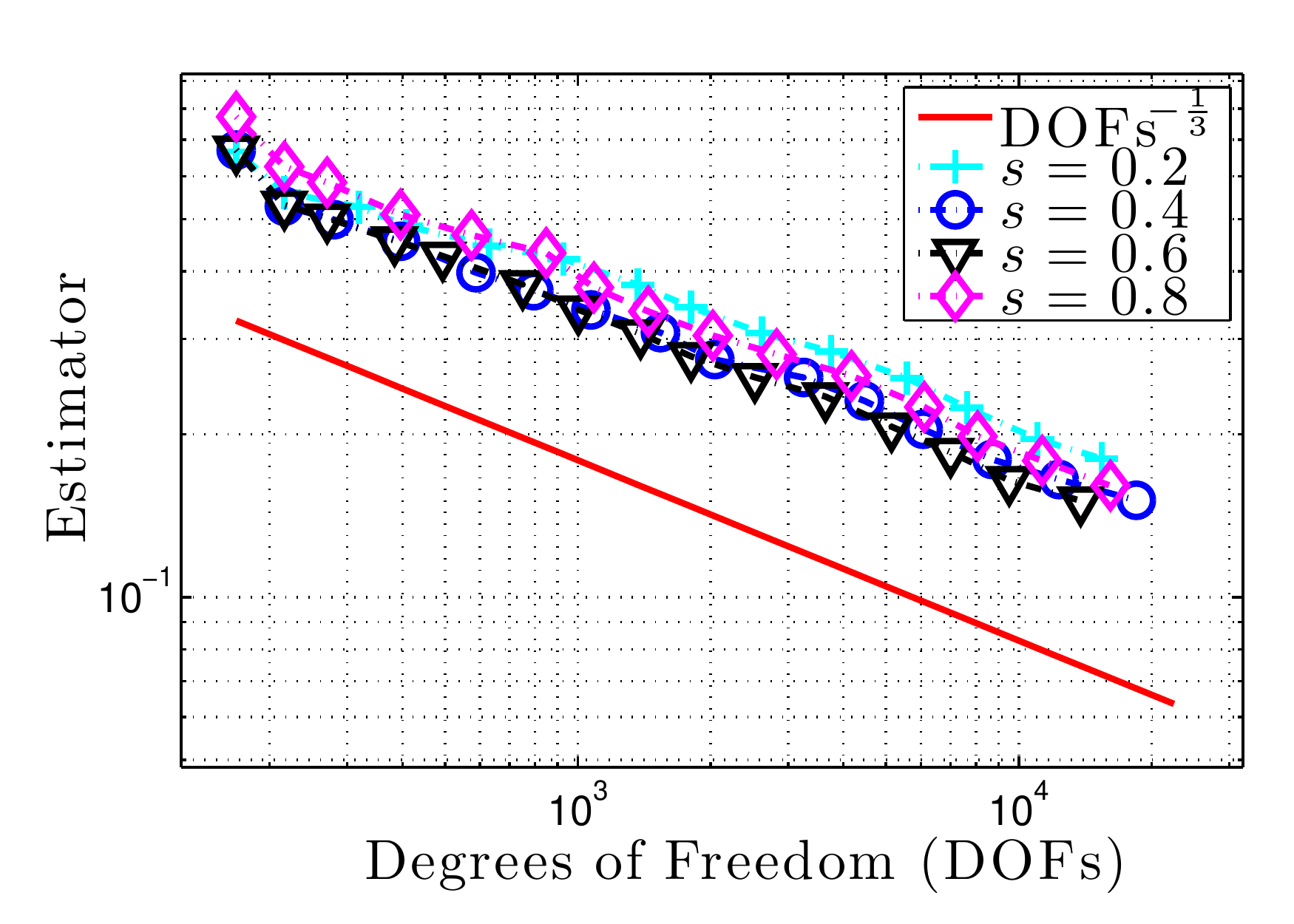}
  \end{center}
  \caption{Computational rate of convergence of AFEM on the smooth but incompatible right hand side over an L-shaped domain of \S~\ref{sub:Lshapedincompatible} for $n=2$ and $s=0.2$, $0.4$, $0.6$ and $s=0.8$. The left panel shows the error vs.~the number of degrees of freedom, whereas the right one that for the total error indicator. In all cases we recover the optimal rate $\# (\T_{\Y})^{-1/3}$.}
\label{fig:Lshapedincompatible}
\end{figure}

We now combine the singularity introduced by the data incompatibility
of \S~\ref{sub:smoothincompatible} and the effect of a reentrant
corner. Consider $\Omega = (-1,1)^2\setminus(0,1)\times(-1,0)$
and $f=1$. As Figure~\ref{fig:Lshapedincompatible} displays, we recover
the  optimal rate of convergence for all possible cases of $s$.

\section{Multilevel methods}
\label{sec:multilevel}

Since we increase the dimension by one in the approximation of
the nonlocal operator $\Laps$ by the local problem \eqref{alpha_harm_intro}, 
it becomes essential to develop efficient linear algebra methods.
It is known that multilevel methods are the most efficient techniques for the solution of discretizations of PDEs \cite{Brandt1977,Brandt1984,Hackbusch1985,Xu1992siamreview}. Multigrid methods for equations of the type \eqref{alpha_harm_intro}, however, have not been explored.

Using the multilevel framework of
\cite{BPWXreg:91,BPWX:91,Xu1992siamreview}, the Xu-Zikatanov identity
\cite{XuZ02} and exploiting the fact that $|y|^\alpha \in A_2(\R^{n+1})$,
which turns out to be critical, we derived in \cite{CNOS} an almost uniform convergent multilevel method to solve \eqref{harmonic_extension_weak}. As we have shown the mesh in the extended dimension must be graded towards the bottom of the cylinder thus becoming anisotropic. We apply line smoothers over vertical lines in the extended domain and prove that the corresponding multigrid \textsf{V}-cycle converges nearly uniformly, \ie the contraction factor depends linearly on the number of levels, and thus logarithmically on the problem size.

\subsection{Multilevel decomposition and multigrid algorithm}
\label{sub:multilevel}

We follow \cite{BP:87,BPWXreg:91} to present a multilevel decomposition of $\V(\T_\Y)$. To simplify the analysis and implementation, we consider a sequence of nested discretizations constructed as follows: We introduce a sequence of nested uniform partitions of the unit interval $\{ \hat{\calI}_k \}$,  with mesh points ${\widehat y}_{l,k}$, for $l=0,\ldots,M_k$ and $k=0,\dots,J$. For $\gamma > 3/(1-\alpha)$ the family $\{\calI_k\}$ is given by
\begin{equation}
\label{graded_mesh_MG}
  y_{l,k} = \Y {\widehat y}_{l,k}^\gamma, \quad l=0,\dots,M_k.
\end{equation}
For $k = 0,\dots,J$, $\T_{\Omega,k}$ is obtained by uniform refinement. The mesh $\T_{\Y,k}$ is the tensor product of $\T_{\Omega,k}$ and $\calI_k$, given by \eqref{graded_mesh_MG}.

If $\V_k:= \V(\T_{\Y,k})$ we have the sequence $\V_0 \subset \V_1 \subset \dots \subset \V_J = \V$, and a space macro-decomposition $\V = \sum_{k=0}^{J} \V_k$. We now introduce a space micro-decomposition: For $j=1,\ldots,\M_k$ let $\calI_{k,j}$ be a subset of $\{1, 2, \ldots , \N_k\}$, and assume that
$
 \cup_{j=1}^{\M_k} \calI_{k,j} = \{1, 2, \ldots , \N_k\}.
$
The sets $\calI_{k,j}$ may not be disjoint but the cardinality of their intersection is bounded independently of $J$ and $\N_{\,k}$. 
If the nodal basis of $\V_k$ is $\phi_{k,i}$, we define $\V_{k,j} = \Span \{ \phi_{k,i}: i\in \calI_{k,j}\}$ and we have the micro-decomposition
\begin{equation}
\label{V=decomp}
  \V = \sum_{k=0}^{J} \sum_{j=1}^{\M_k} \V_{k,j}.
\end{equation} 

With this notation we define a symmetric \textsf{V}-cycle multigrid method as in \cite[Algorithm 3.1]{MR1804746}, with $m \geq 1$ pre and post smoothing steps. When $m=1$, it is equivalent to the application of successive subspace corrections to the decomposition \eqref{V=decomp} with 
exact sub-solvers so that the \textsf{V}-cycle multigrid method has a smoother at each level of block Gauss-Seidel type \cite{BPWXreg:91,Xu1992siamreview}. In particular, the nodal decomposition $\calI_{k,j} = \{j\}$ yields a point-wise Gauss-Seidel smoother. If the indices in $\calI_{k,j}$ are such that the corresponding vertices lie on a vertical line, we obtain \emph{line smoothers}, which are essential to handle anisotropy.

\subsection{Analysis of the multigrid method}
\label{sub:mganalysisdegenerate}

The success of multigrid methods for uniformly elliptic operators is because smoothers are effective in reducing the high frequency components of the error while coarse grid corrections reduce the low frequency ones. However, their effectiveness depends on several factors such as the shape of the mesh. A key ingredient in the design and analysis of a multigrid method on anisotropic meshes is the use of line smoothers \cite{AS:02,BZ:01,Hackbusch:89,S:93}.

When solving \eqref{harmonic_extension_weak} on graded meshes, the approximation on the coarse grid is dominated by the larger mesh size in the $x'$-direction and thus the coarse grid correction cannot capture the smaller scale in the $y$-direction. This is why line smoothing, \ie solving sub-problems restricted to one vertical line is necessary. Owing to the nature of the decomposition, the smoother requires the evaluation of the inverse of the operator over a vertical line. This can be efficiently realized since the corresponding matrix is tridiagonal.

Our analysis hinges on \cite{XuZ02}, which requires stability of the micro-decomposition.

\begin{lemma}[nodal stability and anisotropic inverse inequalities]
\label{lem:anisoInverse}
Let $\T_\Y \in \T$ with $\T_\Omega$ quasiuniform and $\calI_\Y$ graded so that \eqref{graded_mesh_MG} holds. If  $v = \sum_{j=1}^{\M_J} v_j \in \V_J$, then
\begin{equation}
\label{eq:stabaniso}
  \sum_{j=1}^{\M_J} \| v_j \|_{L^2(y^\alpha, \C_\Y)}^2 \lesssim
  \left\| v \right\|_{L^2(y^\alpha, \C_\Y)}^2 \lesssim
  \sum_{j=1}^{\M_J} \| v_j \|_{L^2(y^\alpha, \C_\Y)}^2.
\end{equation}
Moreover, for $T = K \times I \in \T_{\Y}$ we have the following local inverse inequalities
\begin{equation}
\label{eq:invaniso}
  \| \nabla_{x'} v \|_{L^2(y^\alpha,T)} \lesssim h_K^{-1} \| v \|_{L^2(y^\alpha, T)}, \qquad
  \| \partial_y v \|_{L^2(y^\alpha,T)} \lesssim h_I^{-1} \| v \|_{L^2(y^\alpha, T)}.
\end{equation}
\end{lemma}

We now examine the \textsf{V}-cycle multigrid method applied to the decomposition \eqref{V=decomp} with exact sub-solvers on $\V_{k,j}$, \ie with line smoothers; see~\cite[\S{III.12}]{MR1804746} and~\cite{MGaniso}. A key observation in favor of subspaces $\{ \V_{k,j} \}_{j=1}^{\M_k}$ follows.

\begin{lemma}[nodal stability of $y$-derivatives]\label{lem:anisoy}
In the setting of Lemma~\ref{lem:anisoInverse} we have
\begin{equation}
\label{eq:anisoy}
  \sum_{j=1}^{\M_J} \| \partial_y v_j \|_{L^2(y^\alpha, \C_\Y)}^2 \lesssim
  \left\| \partial_y v \right\|_{L^2(y^\alpha, \C_\Y)}^2 \lesssim
  \sum_{j=1}^{\M_J} \| \partial_y v_j \|_{L^2(y^\alpha, \C_\Y)}^2.
\end{equation}
\end{lemma}

We use the quasi-interpolant of \cite{NOS2}, Lemmas~\ref{lem:anisoInverse} and \ref{lem:anisoy} and obtain nearly uniform convergence of the \textsf{V}-cycle multigrid method. We follow \cite{MGaniso,XCN:09}.

\begin{theorem}[convergence of multigrid with line smoothers]
\label{TH:convergence_2}
The symmetric \textup{\textsf{V}}-cycle multigrid method with line smoothing converges with a contraction rate
\[
  \delta \leq 1 - \frac{1}{1+CJ},
\]
where $C$ is independent of $\N_J$ and depends on $y^{\alpha}$ only through $C_{2,y^{\alpha}}$.
\end{theorem}

\subsection{Numerical illustrations}
We consider two cases:
\begin{itemize}
 \item $n=1$, \quad $\Omega = (0,1)$, \quad $u = \sin(3\pi x)$,
 \item $n=2$, \quad $\Omega = (0,1)^2$, \quad $u = \sin (2\pi x_1)\sin(2\pi x_2)$,
\end{itemize}
and $\Y=1$, which has been chosen, as discussed in \cite{NOS}, to capture the exponential decay of the solution.
All of our algorithms are implemented based on the MATLAB$^\copyright$ software package \texttt{{\emph{i}}FEM} \cite{Chen.L2008c}.

Table~\ref{table:1d} shows the number of iterations for the one 
and two dimensional problems, respectively. As we see, the method
converges almost uniformly with respect 
to the number of degrees of freedom. 

\begin{table}[h!]
  \begin{center}
    \begin{tabular}{||l||r||c||c||c||c||}
      \hline
      $h_{\T_{\Omega}}$           & DOFs     & $s=0.15$  & $s=0.3$ & $s=0.6$ & $s=0.8$ \\
      \hline
      $\tfrac1{16}$  & 289     & 7         & 6       & 5       & 5 \\
      \hline
      $\tfrac1{32}$  & 1,089   & 9        & 9       & 6       & 6 \\
      \hline
      $\tfrac1{64}$  & 4,225   & 10        & 10      & 6       & 6 \\
      \hline
      $\tfrac1{128}$ & 16,641  & 10        & 11      & 6       & 6 \\
      \hline
      $\tfrac1{256}$ & 66,049  & 11        & 10      & 6       & 6 \\
      \hline
      $\tfrac1{512}$ & 263,169 & 11        & 10      & 6       & 7 \\
      \hline
    \end{tabular} 
  \end{center}
  \vspace{0.2cm}
  \begin{center}\begin{tabular}{||l||r||c||c||c||c||}
  \hline
  $h_{\T_{\Omega}}$            & DOFs       & $s=0.15$  & $s=0.3$ & $s=0.6$ & $s=0.8$ \\
  \hline
  $\tfrac1{16}$  & 4,913      & 8        & 7       & 6       & 5 \\
  \hline
  $\tfrac1{32}$  & 35,937     & 11        & 8       & 6       & 6 \\
  \hline
  $\tfrac1{64}$  & 274,625    & 12        & 9       & 6       & 6 \\
  \hline
  $\tfrac1{128}$ & 2,146,689  & 13        & 9       & 6       & 6\\
  \hline
  \end{tabular}
  \end{center}
\caption{
Number of iterations for a multigrid method for $(-\Delta)^s u = f$
with $n=1$ (top) and $n=2$ (bottom) using a line smoother 
in the extended direction. The mesh in $\Omega$ is uniform with
mesh size $h_{\T_{\Omega}}$, whereas that in the extended 
direction $(0,\Y)$ is graded according to \eqref{graded_mesh_MG}.}
\label{table:1d}
\end{table}

\section{Space-time fractional diffusion problems}
\label{sec:time_dependent}
We now review the numerical approximation of an initial boundary value problem for a space-time fractional parabolic equation developed in \cite{NOS3}. Given $s \in (0,1)$, $\gamma \in (0,1]$, a function $f$, and an initial datum $u_0$, we seek $u$ such that
\begin{equation}
\label{fractional_heat}
  \partial^{\gamma}_t u + \Laps u = f \ \text{ in } \Omega\times(0,T),
  \quad
  u(0) = u_0 \ \text{ in } \Omega.
\end{equation}
If $\gamma \in (0,1)$, $\partial^{\gamma}_t$ denotes \emph{the left Caputo derivative of order $\gamma$}, defined by
\begin{equation}
\label{caputo}
\partial^{\gamma}_t u(x,t) := \frac{1}{\Gamma(1-\gamma)} \int_{0}^t \frac{1}{(t-r)^{\gamma}} \frac{\partial u(x,r)}{\partial r} \diff r,
\end{equation}
where $\Gamma$ is the Gamma function. For $\gamma = 1$, we consider the usual derivative $\partial_t$.

Fractional derivatives and integrals, like the Caputo derivative, are a powerful tool to describe memory and hereditary properties of materials \cite{MR1926470}.
Problem \eqref{fractional_heat} is motivated by anomalous diffusion processes (subdiffusion) and models dynamical
systems with chaotic motion \cite{MR1604710}, fractional PID controllers \cite{MR1666937}, subdiffusion phenomena in highly heterogeneous aquifers
\cite{PSSB:PSSB2221330150} and plasma turbulence \cite{NEGRETE}.

Problem \eqref{fractional_heat} is complicated due to the nonlocality
of the fractional time derivative \cite{fractional_book,Samko} besides
the presence of the nonlocal space operator. We overcome the latter
with the Caffarelli-Silvestre extension of
\S~\ref{sub:CaffarelliSilvestre} and rewrite \eqref{fractional_heat}
as an elliptic problem with dynamic boundary condition:
\begin{equation}
\label{heat_alpha_extension}
\begin{dcases}
  \DIV \left( y^{\alpha} \nabla \ue \right) = 0 \ \textrm{in } \C\times(0,T), \quad \ue = u_0, \ \textrm{on }\Omega \times \{ 0\}, t=0\\
  \ue = 0 \ \textrm{on }\partial_L \C\times (0,T) \quad d_s \partial_t^{\gamma} \ue + \partial_{\nu}^{\alpha} \ue = d_s f \ \textrm{on } (\Omega \times \{ 0\})\times(0,T).
\end{dcases}
\end{equation}

There are several approaches via finite differences, finite elements and spectral methods to treat the Caputo derivative of order $\gamma$. We refer to \cite[\S1]{MM:13} for an overview. In \cite{LLX:11,LinXu:07} a finite difference scheme is proposed, which has a consistency error $\mathcal{O}(\tau^{2-\gamma})$, where $\tau$ denotes the time step. This estimate, however, requires a rather strong regularity assumption in time \cite{M:10,NOS3}. 
In \cite{NOS3}, we examined the behavior of $\partial_{t} u$ and $\partial_{tt} u$ when $t \downarrow 0$
and derived realistic time-regularity estimates for $u$; see also \cite{M:10,MM:13}. With these refined results we analyzed the truncation error and showed discrete stability. The latter leads to an energy estimate.

\subsection{The Caffarelli-Silvestre extension}
\label{sub:CSparabolic}

We define
\begin{align*}
  \mathcal{W} &:= \{ w \in L^{\infty}(0,T;L^2(\Omega)) \cap L^{2}(0,T;\Hs): \partial_t^{\gamma} 
    w \in L^2(0,T;\Hsd)\}, \\
  \mathcal{V} &:= \{ w \in L^{2}(0,T;\HL(y^{\alpha},\C)): \partial_t^{\gamma} \tr w  \in L^2(0,T;\Hsd)\}.  
\end{align*}
Thus, given $f \in L^2(0,T;\Hsd)$, a function $u \in \mathcal{W}$ solves \eqref{fractional_heat} if and only if the $\alpha$-harmonic extension $\ue \in \mathcal{V}$ solves \eqref{heat_alpha_extension}. A weak formulation of \eqref{heat_alpha_extension} reads: Find $\ue \in \mathcal{V}$ such that $\tr \ue(0) = u_0$ and, for a.e.~$t \in (0,T)$, 
\begin{equation}
\label{heat_harmonic_extension_weak}
\langle \tr  \partial_t^{\gamma} \ue, \tr \phi \rangle +  \frac{1}{d_s} \int_{\C} y^{\alpha} \nabla \ue \nabla \phi = \langle f, \tr \phi \rangle \quad \forall \phi \in \HL(y^{\alpha},\C).
\end{equation}

\begin{remark}[dynamic boundary condition]
\label{rem:dynamical}
Problem \eqref{heat_harmonic_extension_weak} is an elliptic problem with a dynamic boundary condition: $\partial_{\nu}^{\alpha} \ue = d_s f - d_s\tr  \partial_t^{\gamma} \ue$ on $\Omega \times \{0\}$. Its analysis is slightly different from the standard theory for parabolic equations.
\end{remark}

Given $s \in (0,1)$, $\gamma \in (0,1]$, a forcing term $f \in L^2(0,T;\Hsd)$ and $u_0 \in L^2(\Omega)$, problems
\eqref{fractional_heat} and \eqref{heat_alpha_extension} have a unique solution; see \cite[Theorem 2.6]{NOS3}.

\subsection{Time regularity}
\label{sub:Regularity}
As in the elliptic case, in what follows we assume that
\eqref{reg_Omega} holds. We refer to \cite{NOS3} for a complete space
regularity analysis of problem \eqref{fractional_heat}. We focus here on time regularity.

For $\gamma = 1$, we demand sufficient time regularity of the
right-hand side together with compatibility conditions for the initial datum $u_0$. We express this as
\begin{equation}
\label{eq:reg_time_beta01}
  \tr \partial_{tt}\ue \in L^2(0,T;\Hsd).
\end{equation}
For $\gamma \in (0,1)$, \eqref{eq:reg_time_beta01} is inconsistent with the solution of \eqref{fractional_heat}. Properties of the Mittag-Leffler function  show that \eqref{eq:reg_time_beta01} never holds if $u_0 \neq 0$.
Time derivatives of $u$ are unbounded as $t \downarrow 0$. In particular, $\partial_{tt} u (x',t) \notin L^2(0,T;\Hsd)$. However,
\[
  \int_{0^{+}} t^\sigma \|\partial_{tt} u (\cdot,t) \|_{\Hsd}^2 \diff t
\]
is finite provided $\sigma > 3-2\gamma$. For this reason, when $\gamma \in (0,1)$, we assume
\begin{equation*}
t^{\sigma/2} \tr \partial_{tt} \ue \in L^2(0,T;\Hsd) \quad \sigma >
3-2\gamma ,
\end{equation*}
which is a valid assumption provided $\mathcal{A}(u_0,f)<\infty$, where
\begin{equation}
\label{A}
\mathcal{A}(u_0,f) =  \| u_0\|_{\Hs} + \| f\|_{H^2(0,T;\Hsd)}.
\end{equation}
\begin{theorem}[time regularity: $\gamma \in (0,1)$] 
\label{TH:regularity_in_time}
If $f \in H^2(0,T;\Hsd)$ and $u_0 \in \Hs$, then
for $t \in (0,T]$, the solution of \eqref{fractional_heat} satisfies
\begin{equation*}
\| \partial_t u(\cdot,t) - \delta^1 u(\cdot,t) \|_{\Hsd} \lesssim t^{\gamma-1}\mathcal{A}(u_0,f),
\end{equation*}
where $\delta^1 u(\cdot,t) = t^{-1}\big(u(\cdot,t) -u(\cdot,0)\big)$. Moreover,
\begin{equation*}
\| t^{\sigma/2} \partial_{tt} u \|_{L^2(0,T;\Hsd)} \lesssim \mathcal{A}(u_0,f),
\end{equation*}
where $\sigma > 3 -2\gamma$. The hidden constant is independent of $\tau$ but blows up as $\gamma \downarrow0$.
\end{theorem}

\subsection{Time discretization}
\label{sec:time_discretization}

Let $\K \in \mathbb{N}$ denote the number of time steps. The time step is $\tau = T/\K$, and, for $0 \leq k \leq \K$, $t_k = k \tau$ and $I_k= (t_k,t_{k+1}]$. For $\phi \in C( [0,T], \Xcal )$ we denote $\phi^k = \phi(t_k)$ and 
$\phi^{\tau}= \{ \phi^k\}_{k=0}^{\K}$. Moreover, 
\[
\| \phi^{\tau} \|_{\ell^{\infty}(\Xcal)} = \max_{0 \leq k \leq \K} \| \phi^k\|_{\Xcal},
\qquad \| \phi^{\tau} \|_{\ell^2(\Xcal)}^2 = \sum_{k=1}^\K \tau \| \phi^k\|_{\Xcal}^2.
\]
For $W^{\tau} \subset \Xcal$ we define, for $k=0,\dots,\K-1$, $\delta^1 W^{k+1} = \tau^{-1} (W^{k+1} - W^{k})$.

\subsubsection{Time discretization for $\gamma =1$}
\label{sub:discretization_1}
We apply the backward Euler scheme to \eqref{heat_harmonic_extension_weak} for $\gamma=1$: determine $V^\tau = \{ V^k\}_{k=0}^{\K}  \subset \HL(y^{\alpha},\C)$ such that
\begin{equation}
\label{initial_data_cont}
 \tr V^{0} = u_0,
\end{equation}
and, if $f^{k+1} = f(t^{k+1})$ for $k=0,\dots, \K-1$, then $V^{k+1} \in \HL(y^{\alpha},\C)$ solves
\begin{equation*}
   ( \delta^1 \tr  V^{k+1}, \tr W )_{L^2(\Omega)}+
    \frac1{d_s} \int_\C y^\alpha \nabla V^{k+1} \nabla W =  \langle f^{k+1}, \tr W   \rangle, \ \forall W \in \HL(y^{\alpha},\C).
\end{equation*}
Define $U^\tau $ by $U^k= \tr V^k \in \Hs$,
which is a piecewise constant (in time) approximation of $u$, solution to problem \eqref{fractional_heat}. Note that \eqref{initial_data_cont} does not require an extension of $u_0$.
The stability of this scheme is elementary \cite[Lemma 3.3]{NOS3}:
\begin{equation*}
\| \tr V^{\tau} \|^2_{\ell^{\infty}(L^2(\Omega))} + \| V^{\tau} \|^2_{\ell^2(\HLn(y^{\alpha},\C))}
\lesssim  
\| u_0 \|^2_{L^2(\Omega)} + \| f^{\tau} \|^2_{\ell^2(\Hsd)}.
\end{equation*}

\subsubsection{Time discretization for $\gamma \in (0,1)$}
\label{sub:discretization_beta}

We now discretize $\partial_t^{\gamma}$ for $\gamma \in (0,1)$. We consider the 
scheme proposed in \cite{LLX:11,LinXu:07} but resort to the regularity results of Theorem~\ref{TH:regularity_in_time}. Using \eqref{caputo}
and the Taylor formula with integral remainder yields
\begin{equation}
  \begin{aligned}
    \partial_t^{\gamma}u(\cdot,t_{k+1})  =
    \frac{1}{\Gamma(2-\gamma)}\sum_{j=0}^{k} a_j \frac{u(\cdot,t_{k+1-j}) - u(\cdot,t_{k-j})}{ \tau^{\gamma} }
    + \resto_{\gamma}^{k+1}(\cdot),
  \label{discretization_fractional}
  \end{aligned}
\end{equation}
for $0 \leq k \leq \K - 1$, where 
\begin{equation}
\label{a_j}
a_j = (j+1)^{1-\gamma} - j^{1-\gamma}, \qquad 
\resto_{\gamma}^{k+1} = \frac{1}{\Gamma(1-\gamma)} \sum_{j=0}^{k} \int_{I_j} \frac{1}{(t_{k+1}-t)^{\gamma}} R(\cdot,t) \diff t,
\end{equation}
and
\begin{equation}\label{R}
 R(\cdot,t) = \partial_t u(\cdot,t) - \frac{1}{\tau}\big( u(\cdot,t_{j+1}) - u(\cdot,t_j) \big)
\qquad \forall t \in I_j.
\end{equation}
Notice that from \eqref{a_j} we deduce that $a_j > 0$ for all $j\geq 0$ and
$
 1 = a_0 > a_1 > a_2 > \dots > a_j
$, 
$
\lim_{j \rightarrow \infty} a_j = 0
$.

\subsubsection{Consistency}
\label{sub:sub:consistency}
We estimate $\resto_{\gamma}^{\tau}$ using a cancellation property: $R$, defined in \eqref{R}, has vanishing mean in $I_j$ whence we can write
\begin{equation*}
\resto_{\gamma}^{k+1} = \frac{1}{\Gamma(1-\gamma)} \sum_{j=0}^{k} \int_{I_j} ( \psi_{\gamma}(t) - \bar{\psi}^j_{\gamma} )R(\cdot,t) \diff t, 
\end{equation*}
with $\psi_{\gamma}(t) = (t_{k+1}-t)^{-\gamma}$ and
$\bar{\psi}^j_{\gamma} = \fint_{I_j} \psi_{\gamma}(t) \diff
t$. Applying \cite[Lemma 2.1]{NOS3} yields
\begin{equation*}
\| \resto_\gamma^\tau \|_{L^2(0,T;\Hsd)} \lesssim \| \psi_{\gamma} - \bar{\psi}^{\tau}_{\gamma}\|_{L^1(0,T)} \| R^\tau \|_{L^2(0,T;\Hsd)},
\end{equation*}
which reduces the estimation of the residual to deriving suitable bounds for each term on the right hand side of this expression. The regularity results of Theorem \ref{TH:regularity_in_time}, combined with the special structure of the kernel $\psi_{\gamma}$, yield an estimate for $\resto_\gamma^\tau$. This estimate, although giving lower rates of convergence than \cite[(3.4)]{LLX:11}, takes into account the correct behavior of the solution and the singularity of its derivatives as $t\downarrow0$ \cite[Proposition 3.6]{NOS3}:
\begin{equation*}
\| \resto_{\gamma}^{\tau} \|_{L^2(0,T;\Hsd)} \lesssim \tau^{\theta} \mathcal{A}(u_0,f) \quad 0< \theta < 1/2.
\end{equation*}
The hidden constant is independent of the data and $\tau$ but blows up as $\theta \uparrow \frac{1}{2}$.

\subsubsection{Stability and energy estimates}
\label{sub:sub:energy}

We apply \eqref{discretization_fractional}. If $\gamma\in(0,1)$ and $\phi^\tau \subset L^2(\Omega)$ we define the discrete fractional derivative, for $k=0, \ldots, \K-1$, by
\begin{equation*}
\Gamma(2-\gamma)\delta^\gamma \phi^{k+1} := \sum_{j=0}^{k} \frac{  a_j }{ \tau^{\gamma-1} } \delta^1\phi^{k+1-j}
= \frac{\phi^{k+1}}{\tau^\gamma} - \sum_{j=0}^{k-1} \frac{ a_j - a_{j+1} }{\tau^\gamma} \phi^{k-j} - \frac{a_k}{\tau^\gamma} \phi^0.
\end{equation*}
The implicit semi-discrete scheme to solve \eqref{fractional_heat} reads: Set $V^0$ as in \eqref{initial_data_cont} and,
for $k=0,\dots, \K-1$, $V^{k+1} \in \HL(y^\alpha,\C)$ satisfies, for every $W \in \HL(y^\alpha,\C)$,
\begin{equation}
\label{discrete_beta_abs}
(\delta^\gamma \tr V^{k+1}, \tr W)_{L^2(\Omega)} + \frac1{d_s} \int_\C y^\alpha \nabla V^{k+1} \nabla W = \langle f^{k+1},W \rangle.
\end{equation}

Let $I^{\sigma}w$ be the Riemann-Liouville integral of order $\sigma$ of $w$ evaluated at $t=T$.

\begin{theorem}[stability for $\gamma \in (0,1)$]
\label{TH:stabbeta}
The implicit semi-discrete scheme \eqref{discrete_beta_abs} is unconditionally stable and satisfies
\begin{equation}
\label{stab_beta_re}
I^{1-\gamma} \| \tr V^{\tau}\|_{L^2(\Omega)}^2
+ \| V^{\tau} \|^2_{\ell^2(\HLn(y^\alpha,\C))} \leq 
  I^{1-\gamma} \| u_0\|_{L^2(\Omega)}^2 +  \| f^{\tau} \|^2_{\ell^2(\Hsd')}
\end{equation}
\end{theorem}

Deducing an energy estimate for problem \eqref{fractional_heat} is cumbersome due to the nonlocality of the fractional time derivative. A key ingredient in deriving such a result is an integration by parts formula, which for a function not vanishing at $t=0$ and $t=T$ involves boundary terms that need to be estimated; for a step in this direction see \cite{ER:06,LiXu:09}. The energy estimate \eqref{stab_beta_re} yields: 

\begin{corollary}[energy estimate for $u$]
\label{CO:energybeta}
Let $\gamma \in (0,1)$. Then, 
\begin{equation}
\label{stab_beta_re2}
I^{1-\gamma} \|u\|_{L^2(\Omega)}^2
+ \|u \|^2_{L^2(0,T;\Hs)} \leq 
I^{1-\gamma} \| u_0 \|_{L^2(\Omega)}^2 +  \| f \|^2_{L^2(0,T;\Hsd)}.
\end{equation}
\end{corollary}

\begin{remark}[limiting case] 
\label{RE:militing_case}
Given $g \in L^p(0,T)$,  we have $I^\sigma g \rightarrow g$ in $L^p(0,T)$ as $\sigma \downarrow 0$; see \cite[Theorem 2.6]{Samko}. This implies that, taking the limit as $\gamma \uparrow 1$ in \eqref{stab_beta_re2}, we recover the well known stability result for a parabolic equation, \ie
\begin{equation}
\label{stab_beta_1}
\|u\|_{L^\infty(0,T;L^2(\Omega))}^2
+ \|u \|^2_{L^2(0,T;\Hs)} \leq 
\| u_0 \|_{L^2(\Omega)}^2 +  \| f \|^2_{L^2(0,T;\Hsd)}.
\end{equation}
This allows us to unify the estimate of Corollary~\ref{CO:energybeta} for all $\gamma \in (0,1]$.
\end{remark}

\subsection{Truncation}
\label{sec:space_discretization}
We now study the space discretization of \eqref{heat_harmonic_extension_weak}. 
Since $\ue(t)$ decays exponentially in the extended direction $y$, for a.e.~$t \in (0,T)$, we truncate $\C$ to $\C_{\Y}$ for a suitable $\Y$ and seek solutions 
in $\C_{\Y}$.
Define, for $\gamma \in (0,1]$,
\begin{equation}\label{Lambda}
 \Lambda_{\gamma}^2(u_0,f):= I^{1-\gamma}  \| u_0 \|^2_{L^2(\Omega)} 
+ \|f\|^2_{L^2(0,T;\Hsd)},
\end{equation}
where, by  Remark~\ref{RE:militing_case}, $I^0$ is the identity. For $\gamma \in (0,1]$ and $s \in (0,1)$, we have
\begin{equation}
\label{energyYinf}
\|\nabla \ue\|_{L^2\left( 0,T; L^2(y^{\alpha},\Omega \times (\Y,\infty)) \right) } \lesssim e^{-\sqrt{\lambda_1} \Y/2}
\Lambda_{\gamma}(u_0,f),
\end{equation}
where $\Y > 1$ and $\ue$ denotes the solution to \eqref{heat_harmonic_extension_weak}.

As a consequence of \eqref{energyYinf}, we can consider the truncated problem
\begin{equation}
\label{heat_alpha_extension_truncated}
\begin{dcases}
  \DIV \left( y^{\alpha} \nabla v \right) = 0 \ \textrm{in } \C_{\Y}\times(0,T),  
   \quad v = 0 \ \textrm{on } (\partial_L \C_{\Y} \cup \Omega_{\Y}) \times (0,T)
    \\
   d_s \partial_t^{\gamma} \tr v +\partial_{\nu}^{\alpha} v = d_s f \ \textrm{on } (\Omega \times \{ 0\})\times(0,T),
\end{dcases}
\end{equation}
with the initial condition $v = u_0$ on $\Omega \times \{ 0\}$ and $t=0$, and where $\Omega_{\Y} = \Omega \times \{ \Y \}$ with $\Y \geq 1$ sufficiently large. Upon defining
\[
  \mathcal{V}_{\Y}  = \left\{ w \in L^2(0,T;\HL(y^{\alpha},\C_{\Y}) ): \partial_t^{\gamma} \tr w  \in L^2(0,T;\Hsd)\right\},
\]
we understand \eqref{heat_alpha_extension_truncated} as: seek $v \in \mathcal{V}_{\Y}$ such that $\tr v(0) = u_0$ and for a.e.~$t \in (0,T)$,
\begin{equation}
\label{heat_harmonic_extension_weak_truncated}
\langle \partial_t^{\gamma} \tr  v, \tr \phi \rangle + 
\frac1{d_s} \int_{\C_\Y} y^\alpha \nabla v \nabla \phi = \langle f, \tr \phi \rangle,
\quad \forall \phi \in \HL(y^{\alpha},\C_{\Y}).
\end{equation}

If $\ue$ solves \eqref{heat_harmonic_extension_weak}, $v$ solves
\eqref{heat_harmonic_extension_weak_truncated} and
$\Lambda_{\gamma}(u_0,f)$ is given in \eqref{Lambda}, then we have
the following exponential convergence result for every $\gamma \in (0,1]$
and $\Y \geq 1$:
\begin{equation}
\label{exp_convergence}
I^{1-\gamma} \| \tr(\ue-v)\|^2_{L^2(\Omega) }  + 
 \| \nabla(\ue-v)\|^2_{L^2(0,T;L^2(y^\alpha, \C_{\Y}) )} 
\lesssim e^{-\sqrt{\lambda_1} \Y} \Lambda_{\gamma}^2(u_0,f).
\end{equation}

\subsection{Fully discrete scheme}
\label{sec:fully_scheme}

We describe a fully discrete scheme to solve \eqref{heat_harmonic_extension_weak_truncated}. The space discretization hinges on the FEM discussed in Section~\ref{sec:PDE_approach}. The discretization in time uses the schemes proposed in \S\ref{sub:discretization_1} for $\gamma =1$ and in \S\ref{sub:discretization_beta} for $\gamma \in (0,1)$.

The fully discrete scheme computes the sequence $V_{\T_{\Y}}^\tau  \subset \V(\T_{\Y})$, an approximation of the solution to \eqref{heat_harmonic_extension_weak_truncated} at each time step. We first initialize the scheme:
\begin{equation}
\label{initial_data_discrete}
V_{\T_{\Y}}^{0} = \mathcal I_{\T_\Omega}u_0,
\end{equation}
where $\mathcal{I}_{\T_\Omega} = G_{\T_\Y} \circ \mathcal{H}_\alpha$, $\mathcal{H}_\alpha$ denotes the $\alpha$-harmonic extension onto $\C_\Y$
and $G_{\T_\Y}$ the weighted elliptic projection of \cite[\S4.3]{NOS3}; notice that $\tr V_{\T_\Y}^0 = \tr G_{\T_\Y} v(0)$. For $k=0,\dots,\K-1$, let $V_{\T_{\Y}}^{k+1} \in \V(\T_{\Y})$ satisfy for all $W \in \V(\T_\Y)$
 \begin{equation}
 \label{fully_beta}
 ( \delta^\gamma \tr V_{\T_{\Y}}^{k+1} , \tr W )_{L^2(\Omega)}  +
 \frac1{d_s} \int_{\C_\Y} y^\alpha \nabla V_{\T_{\Y}}^{k+1} \nabla W =
 \left\langle f^{k+1}, \tr W   \right\rangle .
\end{equation}
An approximate solution to problem \eqref{fractional_heat} is given by 
$U_{\T_{\Omega}}^\tau = \tr V_{\T_{\Y}}^\tau$.

The discrete scheme \eqref{initial_data_discrete}--\eqref{fully_beta} is unconditionally stable  for all $\gamma \in (0,1]$:
\begin{align}
\nonumber
I^{1-\gamma}\|\tr V_{\T_{\Y}}^\tau \|^2_{L^2(\Omega)}
 & + \| V_{\T_{\Y}}^{\tau} \|^2_{\ell^2(\HLn(y^{\alpha},\C_{\Y}))}
\lesssim
  I^{1-\gamma}\| \tr V^0_{\T_{\Y}} \|_{L^2(\Omega)}^2
+  \| f^{\tau} \|^2_{\ell^2(\Hsd)},
\end{align}
where $I^0$ is the identity according to Remark~\ref{RE:militing_case} (case $\gamma=1$).

We define $\mathcal{B}(\usf_0, f ) := \| \usf_0 \|_{\mathbb{H}^{1+3s}(\Omega)} + \| f|_{t=0} \|_{\mathbb{H}^{1+s}(\Omega)} + \| f \|_{W^1_\infty(0,T;\mathbb{H}^{1-(1-2\mu)s}(\Omega))}$ for $\mu > 0$. The error estimates for \eqref{initial_data_discrete}--\eqref{fully_beta} read as follows.

\begin{theorem}[error estimates for $\gamma \in (0,1)$]
\label{th:order_beta}
Let $\gamma \in (0,1)$, $v$ and $V_{\T_{\Y}}^\tau$ solve \eqref{heat_harmonic_extension_weak_truncated} and \eqref{initial_data_discrete}--\eqref{fully_beta}, respectively. If $\mathcal{A}(u_0,f), \mathcal{B}(u_0,f) < \infty$ and $\T_{\Y}$ satisfies \eqref{graded_mesh}, then
\begin{equation*}
[ I^{1-\gamma}\| \tr (v^\tau - V_{\T_{\Y}}^\tau) \|_{L^2(\Omega)}^2]^{\tfrac{1}{2}} \lesssim 
\tau^{\theta} \mathcal{A}(u_0,f) + |\log N|^{2s} N^\frac{-(1+s)}{n+1} \mathcal{B}(u_0,f),
\end{equation*}
and
\begin{equation*}
\| v^\tau - V_{\T_{\Y}}^\tau \|_{\ell^{2}(\HLn(y^{\alpha},\C_\Y) )} \lesssim \tau^{\theta} \mathcal{A}(u_0,f) 
  + |\log N|^{s} N^\frac{-1}{n+1} \mathcal{B}(u_0,f),
\end{equation*}
where $0<\theta<\frac{1}{2}$, and $\mathcal{A}$ is defined in \eqref{A}.
\end{theorem}

We finally observe that exploiting the regularity \eqref{eq:reg_time_beta01},
we can derive optimal (linear) error estimates in time for $\gamma=1$.

\frenchspacing
\bibliographystyle{plain}
\bibliography{biblio}

\end{document}